\documentclass[a4paper,10pt]{article}
\usepackage{amsmath}
\usepackage{amsfonts}
\usepackage{amssymb}
\usepackage{latexsym}
\usepackage{epsfig}
\usepackage{graphicx}
\usepackage{oldgerm}
\usepackage{theorem}

\def\Z{\mathbb{Z}}
\def\R{\mathbb{R}}
\def\C{\mathbb{C}}
\def\N{\mathbb{N}}
\def\W{\mathbb{W}}
\def\H{\mathbb{H}}

\def\bP{{\bf P}}
\def\bp{{\bf p}}
\def\bE{{\bf E}}

\def\cN{{\cal N}}
\def\cK{{\cal K}}

\def\RN{{\it R}_N}
\def\AN{{\it A}_{N-1}}
\def\BN{{\it B}_N}
\def\BNv{{\it B}^{\vee}_N}
\def\CN{{\it C}_N}
\def\CNv{{\it C}^{\vee}_N}
\def\BCN{{\it BC}_N}
\def\DN{{\it D}_N}

\def\x{\mib{x}}
\def\X{\mib{X}}
\def\v{\mib{v}}
\def\y{\mib{y}}
\def\z{\mib{z}}
\def\bxi{\boldsymbol{\xi}}

\def\rC{{\rm C}}
\def\rc{{\rm c}}
\def\rP{{\rm P}}
\def\rp{{\rm p}}

\def\Det{\mathop{\mathrm{Det}}}

\def\sfp{{\sf p}}
\def\sfM{{\sf M}}
\def\sfr{{\sf r}}
\def\sfA{{\sf A}}

\theorembodyfont{\itshape}

\newtheorem{thm}{Theorem}[section]
\newtheorem{lem}[thm]{Lemma}
\newtheorem{cor}[thm]{Corollary}
\newtheorem{prop}[thm]{Proposition}

\newcommand{\mib}[1]{\mbox{\boldmath $#1$}}
\newcommand{\SSC}[1]{\section{#1}\setcounter{equation}{0}}
\newcommand{\qed}{\hbox{\rule[-2pt]{3pt}{6pt}}}

\begin{document}

\title{\bf 
Macdonald denominators for affine root systems, \\
orthogonal theta functions, \\
and elliptic determinantal point processes
}
\author{
Makoto Katori
\footnote{
Fakult\"{a}t f\"{u}r Mathematik, Universit\"{a}t Wien, 
Oskar-Morgenstern-Platz 1, A-1090 Wien, Austria.
On sabbatical leave from
Department of Physics,
Faculty of Science and Engineering,
Chuo University, 
Kasuga, Bunkyo-ku, Tokyo 112-8551, Japan;
e-mail: katori@phys.chuo-u.ac.jp
}}
\date{28 December 2018}
\pagestyle{plain}
\maketitle

\begin{abstract}
Rosengren and Schlosser introduced notions of
$\RN$-theta functions for the seven types of
irreducible reduced affine root systems,
$\RN=\AN$, $\BN$, $\BNv$, $\CN$, $\CNv$, $\BCN$, $\DN$, $N \in \N$,
and gave the Macdonald denominator formulas.
We prove that, if the variables of the
$\RN$-theta functions are properly scaled
with $N$, they construct seven sets of biorthogonal 
functions, each of which has a continuous parameter
$t \in (0, t_{\ast})$ with given $0< t_{\ast} < \infty$.
Following the standard method in random matrix theory,
we introduce seven types of one-parameter
($t \in (0, t_{\ast})$) families of determinantal point processes
in one dimension,
in which the correlation kernels are expressed by
the biorthogonal theta functions.
We demonstrate that they are elliptic extensions of
the classical determinantal point processes
whose correlation kernels are expressed by
trigonometric and rational functions.
In the scaling limits associated with $N \to \infty$,
we obtain four types of elliptic determinantal point
processes with an infinite number of points
and parameter $t \in (0, t_{\ast})$.
We give new expressions for
the Macdonald denominators using 
the Karlin--McGregor--Lindstr\"om--Gessel--Viennot
determinants for noncolliding Brownian paths,
and show the realization of the associated 
elliptic determinantal point processes
as noncolliding Brownian brides 
with a time duration $t_{\ast}$, which are
specified by the pinned configurations
at time $t=0$ and $t=t_{\ast}$.

\vskip 0.2cm



\end{abstract}
\vspace{3mm}
\normalsize

\SSC
{Introduction} \label{sec:Introduction}

A random $N$-point process, $N \in \N \equiv \{1,2, \dots\}$,
on a space $S \subset \R^d$ 
is a statistical ensemble of nonnegative integer-valued 
Radon measures 
\[
\Xi(\cdot) = \sum_{j=1}^N \delta_{X_j}(\cdot),
\]
where $\delta_y(\cdot), y \in S$ denotes the delta measure such that
$\delta_y(\{x\}) = 1$ if $x=y$ and $\delta_y(\{x\})=0$ otherwise, 
provided that the distribution of points
$\{X_j\}_{j=1}^N$ on $S$ is governed by a probability measure $\bP$.
We assume that $\bP$ has density $\bp$ with respect to
the Lebesgue measure $d \x =\prod_{j=1}^N d x_j$, {\it i.e.}, 
$\bP(\X \in d \x)=\bp(\x) d \x$, $\x \in S^N$. 
For the point process 
$(\Xi, \bP)$,
the {\it $n$-point correlation function}
of a set $\{x_1, \dots, x_n \} \in S^n$, $1 \leq n \leq N$, 
is defined by
\begin{equation}
\rho(\{x_1, \dots, x_n\})
=\frac{1}{(N-n)!} 
\int_{S^{N-n}}
\prod_{j=n+1}^N dx_j \,
\bp(x_1, \dots, x_n, x_{n+1}, \dots, x_N).
\label{eqn:correlation1}
\end{equation}
Then, for any set of observables $\chi_{\ell}, \ell =1,2, \dots, N$,
we have the following useful formulas for expectations, 
\[
\bE \left[
\int_{S^n} \prod_{\ell=1}^n \chi_{\ell}(x_{\ell}) \Xi(dx_{\ell}) \right]
=\int_{S^n} \prod_{\ell=1}^n \Big\{
dx_{\ell} \chi_{\ell}(x_{\ell}) \Big\} 
\rho(\{x_1, \dots, x_n\}),
\quad n =1,2, \dots, N.
\]
If every correlation function is expressed by
a determinant in the form
\begin{equation}
\rho(\{x_1, \dots, x_n\})
=\det_{1 \leq j, k \leq n}[K(x_j, x_k) ]
\label{eqn:correlation2}
\end{equation}
with a two-point continuous function $K(x,y)$, $x, y \in S$,
then the point process is said to be {\it determinantal} 
and $K$ is called the {\it correlation kernel}
\cite{Sos00,ST03,ST03b,AGZ10,Kat15_Springer}. 

A typical example of determinantal point process is the
eigenvalue distribution on $S=\R$ of Hermitian random matrices
in the {\it Gaussian unitary ensemble} (GUE) studied 
in random matrix theory \cite{Meh04,For10,AGZ10}.
The probability measure is given as
\begin{equation}
\bP^{{\rm GUE}_N}(\X \in d \x)=
\bp^{{\rm GUE}_N}(\x) d \x
=\frac{1}{C^{{\rm GUE}_N} }
\prod_{\ell=1}^N e^{-x_{\ell}^2} 
\prod_{1 \leq j < k \leq N} (x_k-x_j)^2 d \x,
\label{eqn:GUE1}
\end{equation}
which is normalized as $(1/N!) \int_{\R^N} \bp^{{\rm GUE}_N}(\x) d \x=1$
with $C^{{\rm GUE}_N}=2^{-N(N-1)/2} \pi^{N/2} \prod_{n=1}^{N-1} n!$.

It is not obvious that one can perform integrations (\ref{eqn:correlation1})
for (\ref{eqn:GUE1}) and obtained results are generally expressed
by determinants as (\ref{eqn:correlation2}).
The verification is, however, not difficult, if we have
the following preliminaries \cite{Meh04,For10,AGZ10,Kat15_Springer}.
\begin{description}
\item{{\bf [P1]}} \,
The factor $\prod_{1 \leq j < k \leq N} (x_k-x_j)$ in (\ref{eqn:GUE1})
obeys the {\it Weyl denominator formula} for the classical root system $\AN$,
\begin{equation}
\det_{1 \leq j, k \leq N} [x_k^{j-1} ] = 
\prod_{1 \leq j < k \leq N} (x_k-x_j).
\label{eqn:Weyl_denominator}
\end{equation}
The determinant in LHS is known as the
{\it Vandermonde determinant}. 
\item{{\bf [P2]}} \,
By a basic property of determinant, without change of value, 
we can replace the entries $x_k^{j-1}$ in LHS
of (\ref{eqn:Weyl_denominator}) 
by any monic polynomials of $x_k$ with order $j-1$.
Here we choose them as the monic Hermitian polynomials
$2^{-(j-1)} H_{j-1}(x) \equiv 2^{-(j-1)} e^{x^2}(-d/d x)^{j-1} e^{-x^2}$, 
and obtain the following equality
including the square roots of Gaussian weights in (\ref{eqn:GUE1}),
\[
\prod_{\ell=1}^N e^{-x_{\ell}^2/2} 
\det_{1 \leq j, k \leq N} [x_k^{j-1} ]
=\det_{1 \leq j, k \leq N} \left[
2^{-(j-1)} e^{-x_k^2/2} H_{j-1}(x_k) \right].
\]
The reason of this choice is that they satisfy
the orthogonality relation, 
\[
\int_{\R} \left\{ 2^{-(j-1)} e^{-x^2/2} H_{j-1}(x) \right\}
\left\{ 2^{-(k-1)} e^{-x^2/2} H_{k-1}(x) \right\}  dx
= h_j \delta_{jk}, \quad
j, k \in \N, 
\]
where
$h_j=2^{-(j-1)} (j-1)! \sqrt{\pi} $.
\end{description}

Then integrals (\ref{eqn:correlation1}) are given by 
determinants (\ref{eqn:correlation2}) with the correlation kernel,
\begin{equation}
K^{{\rm GUE}_N}(x,y)
=\sum_{n=1}^{N} \frac{1}{h_n} 
\left\{ 2^{-(n-1)}e^{-x^2/2} H_{n-1}(x) \right\}
\left\{ 2^{-(n-1)} e^{-y^2/2} H_{n-1}(y) \right\},
\quad x, y \in \R.
\label{eqn:K_GUE}
\end{equation}
See Appendix \ref{sec:appendixC}
for proof in a general setting.

In \cite{RS06}, Rosengren and Schlosser extended 
the Weyl denominator formulas
for classical root systems to the {\it Macdonald denominator formulas}
for seven types of irreducible reduced affine root systems,
$\RN=\AN, \BN, \BNv, \CN, \CNv, \BCN, \DN$ \cite{Mac72}.
They expressed the result using the theta functions
and stated that they are {\it elliptic extensions} of the classical results.
In the present paper, we use their result as
an elliptic extension of the preliminary {\bf [P1]}.
We report in this paper an elliptic extension of the preliminary {\bf [P2]},
and then construct seven types of determinantal 
point processes on the elliptic level, 
$(\Xi^{\RN}, \bP^{\RN}_t, t \in (0, t_{\ast}))$,
$\RN=\AN$, $\BN$, $\BNv$, $\CN$, $\CNv$, $\BCN$, $\DN$,
$N \in \N$,
in the sense that
their correlation kernels are expressed by the orthogonal
theta functions and, if we take appropriate limits of parameters,
they are reduced to the classical ones
expressed by trigonometric and rational functions.

Once the $N$-point systems have been proved to be
determinantal, by taking proper scaling limit
associated with the $N \to \infty$ limit of the correlation kernels,
we can define the determinantal point processes
with an infinite number of points.
Remark that any $N \to \infty$ limit of the probability measure 
$\bP^{{\rm GUE}_N}$ is meaningless, 
since as shown by (\ref{eqn:GUE1})
it is absolutely continuous to the Lebesgue measure of $N$ dimensions,
$d \x=\prod_{j=1}^N dx_j$, 
and $N \to \infty$ limit of $d \x$ cannot be mathematically defined.
Taking the scaling limit of the kernel called the {\it bulk scaling limit},
we obtain the following kernel from (\ref{eqn:K_GUE})
\cite{Meh04,For10,AGZ10,Kat15_Springer}, 
\begin{equation}
K^{\rm sin}(x,y)=\frac{\sin\{ \pi \rho (x-y)\}}{\pi (x-y) },
\quad x, y \in \R.
\label{eqn:sin_kernel0}
\end{equation}
This is called the {\it sine kernel} and it
governs a determinantal point process on $\R$
with an infinite number of points which is
spatially homogeneous on $\R$ with constant density of points $\rho > 0$.

Our elliptic determinantal point processes have
two positive parameters $t_{\ast}$ and $r$.
We demonstrate that in the limit $t_{\ast} \to \infty$,
our seven types of determinantal point processes
on the elliptic level are reduced to the four types
of determinantal point processes on the trigonometric level,
in which the correlation kernels are expressed by
sine functions.
If we take the further limit $r \to \infty$, they are reduced to
the three types of sine kernels, one of which is
identified with (\ref{eqn:sin_kernel0}).
The bulk scaling limit is realized in our systems
by taking the limit
$N \to \infty$, $r \to \infty$ with a ratio $N/r$ fixed
for each $0 < t_{\ast} < \infty$.
We construct four types of 
determinantal point processes on the elliptic level
with an infinite number of particles.
The reductions of them in the limit $t_{\ast} \to \infty$ 
to the classical infinite determinantal point processes are
also shown.

The determinantal point process of GUE, 
$(\Xi^{{\rm GUE}_N}, \bP^{{\rm GUE}_N})$,
is related with an interacting particle system consisting of
$N$ Brownian motions on $\R$, $N \in \N$.
The transition probability density of the one-dimensional
standard Brownian motion (BM) from a point $v$ at time $s$
to a point $x$ at time $t$ is given by
$\rp^{\rm BM}(s, v; t, x)
=e^{-(x-v)^2/\{2(t-s)\}}/\sqrt{2 \pi (t-s)}$,
$0 \leq s < t, v, x \in \R$.
As a function of $x$, 
this is nothing but the probability density function
of the Gaussian distribution with mean $v$ and variance $t-s$.
The square of products of differences,
$\prod_{1 \leq j < k \leq N}(x_k-x_j)^2$, in (\ref{eqn:GUE1})
shows that the points $\{x_j\}_{j=1}^N$ on $\R$
are distributed exclusively.
The corresponding stochastic process is then realized
as a system of Brownian motions conditioned never to collide
with each other \cite{Kat15_Springer}. 
Consider the Weyl chamber,  
$\W_N=\{\x=(x_1, x_2, \dots, x_N) \in \R^N: x_1<x_2< \cdots < x_N\}$.
For $\v, \x \in \W_N$, the total probability mass
of $N$-tuple of Brownian paths, 
in which the $j$-th path starts from $v_j$ at time $s$
and arrives at $x_j$ at time $t >s$, $j=1,2, \dots, N$, 
is given by a determinant $\det[\sfp(s, \v; t, \x) ]$.
Here $\sfp(s, \v; t, \x)$ is the $N \times N$ matrix
whose $(j, k)$-entry is given by
$\rp^{\rm BM}(s, v_j; t, x_k)$,
$j, k \in \{1,2, \dots, N\}$.
This is known as the 
{\it Karlin--McGregor--Lindstr\"om--Gessel--Viennot} (KMLGV)
{\it formula} \cite{KM59,Lin73,GV85}.
Here we consider the situation such that
$N$ BMs start from a given configuration
$\v \in \W_N$ at time 0, execute noncolliding process, and
then return to the configuration $\v$ at time $t_{\ast}>0$.
Such a process is called the $N$-particle system of
{\it noncolliding Brownian bridges} from $\v$ to $\v$ in time duration $t_{\ast}$
(see, for instance, Part I, IV.4.22 of \cite{BS02}
for the original Brownian bridge of a single path). 
The probability density at time $t$ of this $N$-particle process
is then given by (see Section V.C of \cite{KT04})
\[
\bp^{\v \to \v}_t(\x; t_{\ast})
=\frac{\det[\sfp^{\rm BM}(0, \v; t, \x)] \det[\sfp^{\rm BM}(t, \x; t_{\ast}, \v)]}
{\det[\sfp^{\rm BM}(0, \v; t_{\ast}, \v)]},
\quad \x \in \W_N, \quad t \in (0, t_{\ast}).
\]
We can prove that the limit $\v \to {\bf 0} \equiv (0, \dots, 0) \in \R^N$
exists
(see, for instance, Section 3.3 in \cite{Kat15_Springer}),
and we obtain
\begin{equation}
\bp^{{\bf 0} \to {\bf 0}}_t(\x;t_{\ast})
=\frac{1}{C(N, t, t_{\ast})}
\prod_{\ell=1}^N e^{-x_{\ell}^2 t_{\ast}/\{2 t(t_{\ast}-t)\}}
\prod_{1 \leq j < k \leq N} (x_k-x_j)^2,
\quad \x \in \W_N, \quad t \in (0, t_{\ast}), 
\label{eqn:B_bridge2}
\end{equation}
with a normalization factor $C(N, t, t_{\ast})$ which does not depend on $\x$.
If we put $t_{\ast}=2$ and $t=t_{\ast}/2=1$,
(\ref{eqn:B_bridge2}) coincides with 
$\bp^{{\rm GUE}_N}(\x)$ in (\ref{eqn:GUE1}).
In other words, the $N$-particle system
of noncolliding Brownian bridges from ${\bf 0}$ to ${\bf 0}$
with time duration $t_{\ast}$ realizes 
a one-parameter extension 
of determinantal point process of GUE.

Each type of elliptic determinantal point processes
studied in this paper makes a family with one
continuous parameter $t \in (0, t_{\ast})$
(in addition to a discrete parameter $N \in \N$).
We can show that,
$(\Xi^{\AN}, \bP^{\AN}_t, t \in (0, t_{\ast}))$ is realized
as an $N$-particle system of
noncolliding Brownian bridges on 
a circle with radius $r$,
$(\Xi^{\RN}, \bP^{\RN}_t, t \in (0, t_{\ast}))$ 
with $\RN=\BN, \CNv, \BCN$ are
realized as $N$-particle systems of 
noncolliding Brownian bridges in an interval
$[0, \pi r]$ with absorbing boundary condition
at $x=0$ and reflecting boundary condition at $x=\pi r$, 
$(\Xi^{\RN}, \bP^{\RN}_t, t \in (0, t_{\ast}))$ 
with $\RN=\BNv, \CN$ are
realized as $N$-particle systems of 
noncolliding Brownian bridges in an interval
$[0, \pi r]$ with absorbing boundary condition
at both edges, 
and 
$(\Xi^{\DN}, \bP^{\DN}_t, t \in (0, t_{\ast}))$ is
realized as noncolliding $N$-Brownian bridges in 
$[0, \pi r]$ with reflecting boundary condition
at both edges.
These Brownian bridges are specified by the pinned
configurations $\v^{\RN}$ at the initial time $t=0$ and at the final 
time $t=t_{\ast}$.

The paper is organized as follows.
In Section \ref{sec:ortho_theta} we first list out the Macdonald 
denominators $W^{\RN}(\bxi; \tau)$ 
for the seven types of irreducible reduced affine root systems,
$\RN=\AN$, $\BN$, $\BNv$, $\CN$, $\CNv$, $\BCN$, $\DN$, 
and give explicit expressions of theta functions
used by Rosengren and Schlosser 
for the Macdonald denominator formulas \cite{RS06}.
In this paper we use the classical expressions of
Jacobi's theta functions (as shown in Appendix \ref{sec:appendixA})
in order to clarify the conditions that the functions are
real-valued or complex-valued, and to show
dependence on the parameters $t_{\ast}$ and $r$ explicitly.
We prove that, if we put the two variables 
$(\xi, \tau)$ in these theta functions as functions of $(x, t)$
properly depending on the value of $N$,
then the obtained sets of functions
$\{M^{\RN}_j(x, t) \}_{j=1}^N$ constructing 
seven families of biorthogonal systems 
with respect to the integral over $x$
which have a continuous
parameter $t \in (0, t_{\ast})$ (Lemma \ref{thm:orthogonality}).
In Section \ref{sec:DPP} we introduce
seven types of point processes $(\Xi^{\RN}, \bP^{\RN}_t, t \in (0, t_{\ast}))$,
associated with the seven sets of biorthogonal theta functions
after giving the nonnegative conditions (Lemma \ref{thm:positivity})
and the normalization conditions (Lemma \ref{thm:normalization})
for $\bP^{\RN}_t$.
As a byproduct of the latter, the Selberg-type integral 
formulas including the Jacobi theta functions are derived 
as shown in Appendix \ref{sec:appendixB}.
Then we prove that they are all determinantal with
parameter $t \in (0, t_{\ast})$ 
(Theorem \ref{thm:mainA1}). 
The proof of the theorem with derivation of correlation kernels
is given by the standard method in random matrix theory
as explained in Appendix \ref{sec:appendixC}.
We discuss the temporally homogeneous limit
$t_{\ast} \to \infty$ in Section \ref{sec:temp_homo_finite}
and the scaling limit associated with $N \to \infty$ limit
in Section \ref{sec:infinite} (Theorem \ref{thm:infinite_systems}) 
by analyzing the correlation kernels
expressed by $\{M^{\RN}_j(x, t) \}_{j=1}^N$.
Reductions of the determinantal point processes from
the present elliptic level
to the trigonometric and rational function levels
are shown by studying some limit transitions.
In Section \ref{sec:new_exp} 
we give new expressions for the Macdonald denominators
by the KMLGV determinants of
noncolliding BMs (Proposition \ref{thm:KM_det0}).
Then in Section \ref{sec:nB_bridges} 
we show the realizations of $(\Xi^{\RN}, \bP^{\RN}_t, t \in (0, t_{\ast}))$
as $N$-particle systems of Brownian bridges with time 
duration $t_{\ast}$ 
(Theorem \ref{thm:B_bridges}).
Concluding remarks are given in Section \ref{sec:remarks}.

\SSC
{Orthogonal Theta Functions} \label{sec:ortho_theta}
\subsection{Macdonald denominator formulas
of Rosengren and Schlosser}
\label{sec:Macdonald}
Assume that $N \in \N \equiv \{1,2, \dots\}$.
As extensions of the Weyl denominators for
classical root systems, Rosengren and Schlosser \cite{RS06}
studied the Macdonald denominators for the seven
types of irreducible reduced affine root systems \cite{Mac72,Dys72},
$W^{\RN}(\bxi; \tau)$, $\bxi=(\xi_1, \dots, \xi_N) \in \C^N$,
$\RN=\AN$, $\BN, \BNv, \CN, \CNv, \BCN, \DN$,
$N \in \N$.
Up to trivial factors they are written using the Jacobi theta functions
as follows.
(Notations and formulas of the Jacobi theta functions used in this paper
are shown in Appendix \ref{sec:appendixA}.)
\begin{align}
W^{\AN}(\bxi; \tau) &=
\prod_{1 \leq j < k \leq N} \vartheta_1(\xi_k-\xi_j; \tau),
\nonumber\\
W^{\BN}(\bxi; \tau) &=
\prod_{\ell=1}^N \vartheta_1(\xi_{\ell}; \tau)
\prod_{1 \leq j < k \leq N} \Big\{
\vartheta_1(\xi_k-\xi_j; \tau) \vartheta_1(\xi_k+\xi_j; \tau) \Big\},
\nonumber\\
W^{\BNv}(\bxi; \tau) &=
\prod_{\ell=1}^N \vartheta_1(2 \xi_{\ell}; 2 \tau)
\prod_{1 \leq j < k \leq N} \Big\{
\vartheta_1(\xi_k-\xi_j; \tau) \vartheta_1(\xi_k+\xi_j; \tau) \Big\},
\nonumber\\
W^{\CN}(\bxi; \tau) &=
\prod_{\ell=1}^N \vartheta_1(2 \xi_{\ell}; \tau)
\prod_{1 \leq j < k \leq N} \Big\{
\vartheta_1(\xi_k-\xi_j; \tau) \vartheta_1(\xi_k+\xi_j; \tau) \Big\},
\nonumber\\
W^{\CNv}(\bxi; \tau) &=
\prod_{\ell=1}^N \vartheta_1 \left(\xi_{\ell}; \frac{\tau}{2} \right)
\prod_{1 \leq j < k \leq N} \Big\{
\vartheta_1(\xi_k-\xi_j; \tau) \vartheta_1(\xi_k+\xi_j; \tau) \Big\},
\nonumber\\
W^{\BCN}(\bxi; \tau) &=
\prod_{\ell=1}^N \Big\{ \vartheta_1(\xi_{\ell}; \tau) \vartheta_0(2 \xi_{\ell}; 2 \tau) 
\Big\}
\prod_{1 \leq j < k \leq N} \Big\{
\vartheta_1(\xi_k-\xi_j; \tau) \vartheta_1(\xi_k+\xi_j; \tau) \Big\},
\nonumber\\
W^{\DN}(\bxi; \tau) &=
\prod_{1 \leq j < k \leq N} \Big\{
\vartheta_1(\xi_k-\xi_j; \tau) \vartheta_1(\xi_k+\xi_j; \tau) \Big\},
\label{eqn:Macdonald_denominators}
\end{align}
where $\tau \in \H \equiv \{z \in \C : \Im z >0\}$.
They introduced the notions of
$\AN$-theta function of norm $\alpha$ and
$\RN$-theta function for 
$\RN=\BN, \BNv, \CN, \CNv, \BCN, \DN$.
In order to state their results explicitly, here we 
introduce the following four types of functions,
\begin{align}
\Theta^{A}(\sigma, z, \tau) 
&=e^{2 \pi i \sigma z} 
\vartheta_2 (\sigma \tau + z; \tau),
\nonumber\\
\Theta^{B}(\sigma, z, \tau) 
&= e^{2 \pi i \sigma z} \vartheta_1 (\sigma \tau + z; \tau) 
- e^{-2 \pi i \sigma z} \vartheta_1 (\sigma \tau - z ; \tau), 
\nonumber\\
\Theta^{C}(\sigma, z, \tau) 
&= e^{2 \pi i \sigma z} \vartheta_2 (\sigma \tau + z; \tau)
- e^{-2 \pi i \sigma z} \vartheta_2 (\sigma \tau - z; \tau),
\nonumber\\
\Theta^{D}(\sigma, z, \tau)  
&= e^{2 \pi i \sigma z} \vartheta_2 (\sigma \tau + z; \tau)
+ e^{-2 \pi i \sigma z} \vartheta_2 (\sigma \tau - z; \tau),
\label{eqn:Theta}
\end{align}
for $\sigma \in \R, z \in \C, \tau \in \H$,
where $i=\sqrt{-1}$. 
Let
\[
\sharp(\RN)
= \begin{cases}
A, \quad & \mbox{if $\RN=\AN$},
\cr
B, \quad & \mbox{if $\RN=\BN, \BNv$},
\cr
C, \quad & \mbox{if $\RN=\CN, \CNv, \BCN$},
\cr
D, \quad & \mbox{if $\RN=\DN$},
\end{cases}
\]
\begin{equation}
J^{\RN}(j) = \begin{cases}
j-1/2, \quad & \RN = \AN, \CNv,
\\
j-1, \quad & \RN=\BN, \BNv, \DN,
\\
j, \quad & \RN=\CN, \BCN,
\end{cases}
\label{eqn:J_R}
\end{equation}
and
\begin{equation}
\cN^{\RN} = \begin{cases}
N, \quad & \RN = \AN, \\
2N-1, \quad & \RN = \BN, \\
2N, \quad & \RN = \BNv, \CNv, \\
2(N+1), \quad & \RN =\CN, \\
2N+1, \quad & \RN = \BCN, \\
2(N-1), \quad & \RN = \DN.
\end{cases}
\label{eqn:N_R}
\end{equation}
Rosengren and Schlosser proved that, if we put
\begin{equation}
f^{\RN}_j(\xi; \tau)
=\Theta^{\sharp(\RN)}(J^{\RN}(j), \xi, \cN^{\RN} \tau),
\quad j=1,2, \dots, N,
\label{eqn:RS_AAA}
\end{equation}
then 
\begin{equation}
\det_{1 \leq j, k \leq N}
\left[ f^{\AN}_{j}(\xi_k; \tau) \right]
= C^{\AN}(\tau) \vartheta_1 \left( \sum_{\ell=1}^N \xi_{\ell}+\widetilde{\alpha} \right)
W^{\AN}(\bxi; \tau)
\label{eqn:Macdonald1}
\end{equation}
with norm
$\alpha=e^{2 \pi i \widetilde{\alpha}}$,
and 
\begin{equation}
\det_{1 \leq j, k \leq N}
\left[ f^{\RN}_{j}(\xi_k; \tau) \right]
= C^{\RN}(\tau) W^{\RN}(\bxi; \tau),
\quad \mbox{for $\RN=\BN, \BNv, \CN, \CNv, \BCN, \DN$},
\label{eqn:Macdonald2}
\end{equation}
where $C^{\RN}(\tau)$ depend on $\tau$, $N$
(and the choice of norm $\alpha$ for $\AN$), 
but not on $\bxi$.
The factors $C^{\RN}(\tau)$ are explicitly given
in Proposition 6.1 in \cite{RS06} 
and the equalities (\ref{eqn:Macdonald1}) and (\ref{eqn:Macdonald2})
are called the {\it Macdonald denominator formulas}.
See also \cite{Kra05,War02}.

\subsection{Biorthogonality}
\label{sec:orthogonality}

Assume that $0< r < \infty$.
Let 
\begin{equation}
\xi(x) = \xi(x;r)= \frac{x}{2 \pi r},
\qquad
\tau(t)=\tau(t; r)= \frac{i t}{2 \pi r^2}. 
\label{eqn:xi_tau1}
\end{equation}
In the present paper, we consider the following
seven sets of functions of
$(x, t) \in \R \times [0, \infty)$,
$\{M^{\RN}_j(x,t) \}_{j=1}^N$,
which are defined using the $\AN$-theta function of
norm $\alpha=e^{2 \pi i \widetilde{\alpha}_N}$ with
\begin{equation}
\widetilde{\alpha}_N= \begin{cases}
N \tau(t)/2, & \mbox{if $N$ is even},
\cr
(1+N \tau(t))/2, & \mbox{if $N$ is odd},
\end{cases}
\label{eqn:norm}
\end{equation}
and the $\RN$-theta functions,
$\RN=\BN, \BNv, \CN, \CNv, \BCN, \DN$,
of Rosengren and Schlosser as
\begin{equation}
M^{\RN}_j(x, t)
= \Theta^{\sharp(\RN)} 
\Big(J^{\RN}(j)/\cN^{\RN}, \cN^{\RN} \xi(x); (\cN^{\RN})^2 \tau(t) \Big), \quad
j=1,2, \dots, N.
\label{eqn:MR0}
\end{equation}
Note that the choice of norm (\ref{eqn:norm}) for the
$\AN$-theta function is different from the
previous papers \cite{Kat15,Kat16,Kat17}.

The explicit expressions of these functions are
given by follows, 
\begin{align}
& M^{\AN}_j(x, t) = M^{\AN}_j(x, t; r)
\nonumber\\
& \qquad 
=e^{2 \pi i J^{\AN}(j) \xi(x)} 
\vartheta_2 \Big( \cN^{\AN} \{J^{\AN}(j) \tau(t) + \xi(x)\}; (\cN^{\AN})^2 \tau(t) \Big),
\label{eqn:MA1}
\\
& M^{\RN}_j(x, t) = M^{\RN}_j(x, t; r)
\nonumber\\
& \qquad 
= e^{2 \pi i J^{\RN}(j) \xi(x)} \vartheta_1 \Big(
\cN^{\RN} \{J^{\RN}(j) \tau(t) + \xi(x)\}; (\cN^{\RN})^2 \tau(t) \Big)
\nonumber\\
& \qquad \quad 
- e^{-2 \pi i J^{\RN}(j) \xi(x)} \vartheta_1 \Big(
\cN^{\RN} \{J^{\RN}(j) \tau(t) - \xi(x)\}; (\cN^{\RN})^2 \tau(t) \Big), 
\quad 
\mbox{for $\RN=\BN, \BNv$},
\label{eqn:MB1}
\\
& M^{\RN}_j(x, t) = M^{\RN}_j(x, t; r)
\nonumber\\
& \qquad
= e^{2 \pi i J^{\RN}(j) \xi(x)} \vartheta_2 \Big(
\cN^{\RN} \{J^{\RN}(j) \tau(t) + \xi(x)\}; (\cN^{\RN})^2 \tau(t) \Big)
\nonumber\\
& \qquad \quad 
- e^{-2 \pi i J^{\RN}(j) \xi(x)} \vartheta_2 \Big(
\cN^{\RN} \{J^{\RN}(j) \tau(t) - \xi(x)\}; (\cN^{\RN})^2 \tau(t) \Big),
\quad 
\mbox{for $\RN=\CN, \CNv, \BCN$},
\label{eqn:MC1}
\\
&M^{\DN}_j(x, t) = M^{\DN}_j(x, t; r)
\nonumber\\
& \qquad 
= e^{2 \pi i J^{\DN}(j) \xi(x)} \vartheta_2 \Big(
\cN^{\DN} \{J^{\DN}(j) \tau(t) + \xi(x)\}; (\cN^{\DN})^2 \tau(t) \Big)
\nonumber\\
& \qquad \quad
+ e^{-2 \pi i J^{\DN}(j) \xi(x)} \vartheta_2 \Big( 
 \cN^{\DN} \{J^{\DN}(j)\tau(t) - \xi(x)\}; (\cN^{\DN})^2 \tau(t) \Big).
\label{eqn:MD1}
\end{align}
The complex conjugates of these functions are given as
\begin{align}
\overline{M^{\AN}_j(x,t)}
&=M^{\AN}_j(-x, t), 
\nonumber\\
\overline{M^{\RN}_j(x,t)}
&= M^{\RN}_j(x, t) \in \R, \quad \mbox{for $\RN=\BN,\BNv,\DN$},
\nonumber\\
\overline{M^{\RN}_j(x,t)}
&= - M^{\RN}_j(x, t) \in i \R, \quad \mbox{for $\RN=\CN,\CNv,\BCN$}.
\label{eqn:complex_conjugates}
\end{align}

The setting of three variables of
$\Theta^{\sharp(\RN)}$ in (\ref{eqn:MR0}), 
which is different from that in 
(\ref{eqn:RS_AAA}), is essential
for establishing the following biorthogonality relations.

\begin{lem}
\label{thm:orthogonality}
Assume $0 < t_{\ast} < \infty$.
For any $t \in (0, t_{\ast})$, if $j, k \in \{1,2, \dots, N\}$, then
\begin{align}
& \int_0^{2 \pi r} \overline{M^{\AN}_j(x, t_{\ast}-t)}
M^{\AN}_k(x, t) dx
= m^{\AN}_j(t_{\ast}) \delta_{jk},
\label{eqn:ortho_A1}
\\
& \int_0^{\pi r} \overline{M^{\RN}_j(x, t_{\ast}-t)}
M^{\RN}_k(x, t) dx
= m^{\RN}_j(t_{\ast}) \delta_{jk},
\quad \mbox{for $\RN=\BN,\BNv,\CN,\CNv,\BCN,\DN$}, 
\label{eqn:ortho_BD1}
\end{align}
where
\begin{align}
m^{\RN}_j(t_{\ast}) &=
2 \pi r \vartheta_2 \Big( \cN^{\RN} J^{\RN}(j) \tau(t_{\ast}); (\cN^{\RN})^2 \tau(t_{\ast}) \Big),
\quad j \in \{1,2, \dots, N \},
\nonumber\\
& \hskip 8.5cm
\mbox{for $\RN=\AN, \CN, \CNv, \BCN$},
\label{eqn:ortho_const1}
\\
m^{\RN}_j(t_{\ast}) &= \begin{cases}
4 \pi r \vartheta_2 \Big(0; (\cN^{\RN})^2 \tau(t_{\ast}) \Big),
& j=1, \\
2 \pi r \vartheta_2 \Big(\cN^{\RN} J^{\RN}(j) \tau(t_{\ast});
 (\cN^{\RN})^2 \tau(t_{\ast}) \Big),
& j \in \{2, 3, \dots, N\},
\end{cases}
\quad \mbox{for $\RN=\BN, \BNv$},
\label{eqn:ortho_const2}
\\
m^{\DN}_j(t_{\ast}) &= \begin{cases}
4 \pi r \vartheta_2 \Big(0; (\cN^{\DN})^2 \tau(t_{\ast}) \Big),
& j=1,
\\
2 \pi r \vartheta_2 \Big(\cN^{\DN} J^{\DN}(j) \tau(t_{\ast});
 (\cN^{\DN})^2 \tau(t_{\ast}) \Big),
& j \in \{2, 3, \dots, N-1\},
\\
4 \pi r \vartheta_2 \Big(\cN^{\DN}(N-1) \tau(t_{\ast}); 
(\cN^{\DN})^2 \tau(t_{\ast}) \Big),
& j=N.
\end{cases}
\label{eqn:ortho_const3}
\end{align}
\end{lem}
{\it Proof.} \,
(i) First we prove (\ref{eqn:ortho_A1}) for the type $\AN$.
By (\ref{eqn:J_R}), (\ref{eqn:N_R}), and (\ref{eqn:complex_conjugates}), 
\begin{align*}
I^{\AN}_{jk} \equiv& \int_0^{2 \pi r} \overline{M^{\AN}_j(x, t_{\ast}-t)}
M^{\AN}_k(x, t) dx
\nonumber\\
&= \int_{0}^{2 \pi r} dx \,
e^{2 \pi i(k-j) \xi(x)}
\vartheta_2 \left( N \{(j-1/2) \tau(t_{\ast}-t) - \xi(x) \};
N^2 \tau(t_{\ast}-t) \right)
\nonumber\\
& \qquad \qquad \qquad \times
\vartheta_2\left( N \{ (k-1/2) \tau(t)+ \xi(x) \};
N^2 \tau(t) \right)
\nonumber\\
&= \sum_{n \in \Z} \sum_{m \in \Z}
e^{
(n-1/2)^2 N^2 \tau(t_{\ast}-t) \pi i + (2n-1) N (j-1/2) \tau(t_{\ast}-t) \pi i}
\nonumber\\
& \qquad \qquad \times
e^{
(m-1/2)^2 N^2 \tau(t) \pi i + (2m-1)N  (k-1/2) \tau(t) \pi i }
\int_0^{2 \pi r} e^{2 \pi i \{(k-j)- N(n-m)\}\xi(x)} dx,
\end{align*}
where we have used
the definition of $\vartheta_2$ given by (\ref{eqn:theta}).
By (\ref{eqn:xi_tau1}), 
\[
\int_0^{2 \pi r} e^{2 \pi i \{(k-j)-N(n-m)\} \xi(x)} dx
= 2 \pi r \int_0^1 e^{2 \pi i \{(k-j)-N(n-m)\} \xi} d \xi.
\]
Here we use the equality
\[
\int_0^{1} e^{ 2 \pi i \theta \xi} d \xi = {\bf 1}(\theta=0),
\]
where ${\bf 1}(\omega)$ is 
the indicator function of condition $\omega$;
${\bf 1}(\omega)=1$ if $\omega$ is satisfied,
and ${\bf 1}(\omega)=0$ otherwise. 
The integral $I_{jk}^{\AN}$ is nonzero, if and only if
\begin{equation}
(k-j)- N (n-m) =0.
\label{eqn:cond_A1}
\end{equation}
Since $j, k \in \{1, 2, \dots, N\}$ and $n, m \in \Z$,
we see that
$-(N-1) \leq k-j \leq N-1$, and 
$N (n-m) \in N \Z$.
Hence (\ref{eqn:cond_A1}) is satisfied if and only if,
$j=k$ and $n=m$.
Therefore, we can conclude $I^{\AN}_{jk}=0$, if $j \not=k$,
and
\begin{align*}
I^{\AN}_{jj} &= 2 \pi r \sum_{n \in \Z} 
e^{\{ (n-1/2)^2 N^2 
+(2n-1) N(j-1/2)\} \{\tau(t_{\ast}-t)+\tau(t) \} \pi i }
\nonumber\\
&= 2 \pi r \vartheta_2( N (j-1/2) \tau(t_{\ast}); N^2 \tau(t_{\ast})).
\end{align*}
Then for $\RN=\AN$ the proof of (\ref{eqn:ortho_A1}) 
with (\ref{eqn:ortho_const1})
is complete. \\
(ii) Next we prove (\ref{eqn:ortho_BD1}) for $\RN=\BN$ and $\BNv$.
By (\ref{eqn:complex_conjugates}), LHS of (\ref{eqn:ortho_BD1})
is given by
\begin{align*}
I^{\RN}_{jk} &= \int_0^{\pi r} M^{\RN}_j(x, t_{\ast}-t) M^{\RN}_k(x, t) dx
\nonumber\\
&= I^{\RN}_{jk, ++}-I^{\RN}_{jk, +-} - I^{\RN}_{jk, -+} + I^{\RN}_{jk, --},
\end{align*}
where
\begin{align*}
I^{\RN}_{jk, \pm \pm}
&= \int_0^{\pi r} dx \, e^{2 \pi i \{\pm(j-1)\pm(k-1)\} \xi(x)}
\vartheta_1 \Big(\cN^{\RN}\{(j-1) \tau(t_{\ast}-t) \pm \xi(x) \};
(\cN^{\RN})^2 \tau(t_{\ast}-t) \Big)
\nonumber\\
& \qquad \qquad \qquad \times
\vartheta_1 \Big( \cN^{\RN}\{(k-1) \tau(t) \pm \xi(x) \};
(\cN^{\RN})^2 \tau(t) \Big).
\end{align*}
By changing the sign of integral variables appropriately, $x \to - x$, we obtain
\begin{align*}
I^{\RN}_{jk, +} 
&\equiv I^{\RN}_{jk, ++}+I^{\RN}_{jk, --}
\nonumber\\
&= \int_{-\pi r}^{\pi r} dx \, e^{2 \pi i (j+k-2) \xi(x)}
\vartheta_1 \Big( \cN^{\RN}\{(j-1) \tau(t_{\ast}-t) + \xi(x) \};
(\cN^{\RN})^2 \tau(t_{\ast}-t) \Big)
\nonumber\\
& \qquad \qquad \qquad \times
\vartheta_1 \Big( \cN^{\RN}\{(k-1) \tau(t) + \xi(x) \};
(\cN^{\RN})^2 \tau(t) \Big),
\nonumber\\
I^{\RN}_{jk, -} 
&\equiv I^{\RN}_{jk, +-}+I^{\RN}_{jk, -+}
\nonumber\\
&= \int_{-\pi r}^{\pi r} dx \, e^{2 \pi i (j-k) \xi(x)}
\vartheta_1 \Big( \cN^{\RN}\{(j-1) \tau(t_{\ast}-t) + \xi(x) \};
(\cN^{\RN})^2 \tau(t_{\ast}-t) \Big)
\nonumber\\
& \qquad \qquad \qquad \times
\vartheta_1 \Big( \cN^{\RN}\{(k-1) \tau(t) - \xi(x) \};
(\cN^{\RN})^2 \tau(t) \Big).
\end{align*}
By the definition of $\vartheta_1$ given by (\ref{eqn:theta}), we have
\begin{align*}
I^{\RN}_{jk, +}
&= - \sum_{n \in \Z} \sum_{m \in \Z}
(-1)^{n+m} e^{ (n-1/2)^2 (\cN^{\RN})^2 \tau(t_{\ast}-t) \pi i
+(2n-1) \cN^{\RN}(j-1) \tau(t_{\ast}-t) \pi i }
\nonumber\\
& \qquad \times
e^{
(m-1/2)^2 (\cN^{\RN})^2 \tau(t) \pi i + (2m-1) \cN^{\RN} (k-1) \tau(t) \pi i }
\int_{-\pi r}^{\pi r} e^{2 \pi i \{(j+k-2)+\cN^{\RN}(n+m-1)\}\xi(x)} dx,
\nonumber\\
I^{\RN}_{jk, -}
&= - \sum_{n \in \Z} \sum_{m \in \Z}
(-1)^{n+m} e^{(n-1/2)^2 (\cN^{\RN})^2 \tau(t_{\ast}-t) \pi i
+(2n-1) \cN^{\RN}(j-1) \tau(t_{\ast}-t) \pi i }
\nonumber\\
& \qquad \times
e^{
(m-1/2)^2 (\cN^{\RN})^2 \tau(t) \pi i + (2m-1) \cN^{\RN} (k-1) \tau(t) \pi i }
\int_{-\pi r}^{\pi r} e^{2 \pi i \{(j-k)+\cN^{\RN}(n-m)\}\xi(x)} dx.
\end{align*}
Here we note that
\begin{align*}
\int_{-\pi r}^{\pi r} e^{2 \pi i \{(j+k-2)+\cN^{\RN}(n+m-1)\} \xi(x)} dx
&= 2 \pi r \int_{-1/2}^{1/2} e^{2 \pi i \{(j+k-2)+\cN^{\RN}(n+m-1)\} \xi} d \xi
\nonumber\\
&= 2 \pi r {\bf 1}\Big( (j+k-2)+\cN^{\RN}(n+m-1)=0 \Big),
\nonumber\\
\int_{-\pi r}^{\pi r} e^{2 \pi i \{(j-k)+\cN^{\RN}(n-m)\} \xi(x)} dx
&= 2 \pi r \int_{-1/2}^{1/2} e^{2 \pi i \{(j-k)+\cN^{\RN}(n-m)\} \xi} d \xi
\nonumber\\
&= 2 \pi r {\bf 1}\Big( (j-k)+\cN^{\RN}(n-m)=0 \Big).
\end{align*}
Since $j, k \in \{1,2, \dots, N\}$, we see that
\[
0 \leq j+k-2 \leq 2N-2 < \cN^{\RN}
=\begin{cases} 
2N-1, & \RN=\BN, \cr
2N, & \RN=\BNv,
\end{cases}
\]
and
$0 \leq j-k \leq N-1 < \cN^{\RN}$.
The condition $(j+k-2)+\cN^{\RN}(n+m-1)=0$ is satisfied if and only if
\[
j+k-2=0, \quad n+m-1=0
\quad \Longleftrightarrow \quad
j=k=1, \quad m=-n+1,
\]
and the condition $(j-k)+\cN^{\RN}(n-m)=0$ is
satisfied if and only if
$j=k, n=m$.
Hence
\begin{align*}
I^{\RN}_{jk,+}
&= - 2 \pi r \delta_{j1} \delta_{k1}
\sum_{n \in \Z} (-1)
e^{(n-1/2)^2 (\cN^{\RN})^2 \tau(t_{\ast}-t) \pi i }
\times
e^{ (-n+1/2)^2 (\cN^{\RN})^2 \tau(t) \pi i }
\nonumber\\
&= 2 \pi r \delta_{jk} \delta_{j1} \vartheta_2 \Big( 0; (\cN^{\RN})^2 \tau(t_{\ast}) \Big),
\end{align*}
and
$I^{\RN}_{jk,-}
= - 2 \pi r \delta_{jk} \vartheta_2 ( \cN^{\RN}(j-1) \tau(t_{\ast}); (\cN^{\RN})^2 \tau(t_{\ast}) )$.
Therefore, we obtain
\begin{align*}
I^{\RN}_{jk} &= I^{\RN}_{jk, +}-I^{\RN}_{jk, -}
\nonumber\\
&= 2 \pi r \delta_{jk}
\Big\{ \delta_{j1}
\vartheta_2 \Big( 0; (\cN^{\RN})^2 \tau(t_{\ast}) \Big)
+
\vartheta_2 \Big( \cN^{\RN}(j-1) \tau(t_{\ast}); (\cN^{\RN})^2 \tau(t_{\ast}) \Big)
\Big\}.
\end{align*}
This proves (\ref{eqn:ortho_BD1}) for $\RN=\BN$ and $\BNv$
with (\ref{eqn:ortho_const2}).
\\
(iii) Now we prove (\ref{eqn:ortho_BD1}) for $\RN=\CN$ and $\BCN$.
By (\ref{eqn:complex_conjugates}), LHS of (\ref{eqn:ortho_BD1})
is given by
\[
I^{\RN}_{jk} = - \int_0^{\pi r} M^{\RN}_j(x, t_{\ast}-t) M^{\RN}_k(x, t) dx
= - I^{\RN}_{jk, +} + I^{\RN}_{jk, -},
\]
where
\begin{align*}
I^{\RN}_{jk, +} 
&= \int_{-\pi r}^{\pi r} dx \, e^{2 \pi i (j+k) \xi(x)}
\vartheta_2 \Big( \cN^{\RN}\{j \tau(t_{\ast}-t) + \xi(x) \};
(\cN^{\RN})^2 \tau(t_{\ast}-t) \Big)
\nonumber\\
& \qquad \qquad \qquad \times
\vartheta_2 \Big( \cN^{\RN}\{k \tau(t) + \xi(x) \};
(\cN^{\RN})^2 \tau(t) \Big),
\nonumber\\
I^{\RN}_{jk, -} 
&= \int_{-\pi r}^{\pi r} dx \, e^{2 \pi i (j-k) \xi(x)}
\vartheta_2 \Big( \cN^{\RN}\{j \tau(t_{\ast}-t) + \xi(x) \};
(\cN^{\RN})^2 \tau(t_{\ast}-t) \Big)
\nonumber\\
& \qquad \qquad \qquad \times
\vartheta_2 \Big( \cN^{\RN}\{k \tau(t) - \xi(x) \};
(\cN^{\RN})^2 \tau(t) \Big).
\end{align*}
By the definition of $\vartheta_2$ given by (\ref{eqn:theta}), we have
\begin{align*} 
I^{\RN}_{jk, +}
&= \sum_{n \in \Z} \sum_{m \in \Z}
e^{ (n-1/2)^2 (\cN^{\RN})^2 \tau(t_{\ast}-t) \pi i
+(2n-1) \cN^{\RN} j \tau(t_{\ast}-t) \pi i }
\nonumber\\
& \qquad \times
e^{
(m-1/2)^2 (\cN^{\RN})^2 \tau(t) \pi i + (2m-1) \cN^{\RN} k \tau(t) \pi i }
\int_{-\pi r}^{\pi r} e^{2 \pi i \{(j+k)+\cN^{\RN}(n+m-1)\}\xi(x)} dx,
\nonumber\\
I^{\RN}_{jk, -}
&= \sum_{n \in \Z} \sum_{m \in \Z}
e^{(n-1/2)^2 (\cN^{\RN})^2 \tau(t_{\ast}-t) \pi i
+(2n-1) \cN^{\RN}j \tau(t_{\ast}-t) \pi i }
\nonumber\\
& \qquad \times
e^{
(m-1/2)^2 (\cN^{\RN})^2 \tau(t) \pi i + (2m-1) \cN^{\RN} k \tau(t) \pi i }
\int_{-\pi r}^{\pi r} e^{2 \pi i \{(j-k)+\cN^{\RN}(n-m)\}\xi(x)} dx.
\end{align*}
Here we note that
\begin{align*}
\int_{-\pi r}^{\pi r} e^{2 \pi i \{(j+k)+\cN^{\RN}(n+m-1)\} \xi(x)} dx
&= 2 \pi r {\bf 1}\Big( (j+k)+\cN^{\RN}(n+m-1)=0 \Big),
\nonumber\\
\int_{-\pi r}^{\pi r} e^{2 \pi i \{(j-k)+\cN^{\RN}(n-m)\} \xi(x)} dx
&= 2 \pi r {\bf 1}\Big( (j-k)+\cN^{\RN}(n-m)=0 \Big).
\end{align*}
Since $j, k \in \{1,2, \dots, N\}$, we see that
\[
2 \leq j+k \leq 2N < \cN^{\RN}
=\begin{cases} 
2(N+1), & \RN=\CN, \cr
2N+1, & \RN=\BCN,
\end{cases}
\]
and
$0 \leq j-k \leq N-1 < \cN^{\RN}$.
The condition $(j+k)+\cN^{\RN}(n+m-1)=0$ is not satisfied,
and thus $I^{\RN}_{jk,+} \equiv 0$.
The condition $(j-k)+\cN^{\RN}(n-m)=0$ is
satisfied if and only if $j=k, n=m$.
Hence
\begin{align*}
I^{\RN}_{jk} &= I^{\RN}_{jk, -}
=2 \pi r \delta_{jk}
\sum_{n \in \Z} 
e^{ (n-1/2)^2 (\cN^{\RN})^2 \tau(t_{\ast}) \pi i 
+(2n-1) \cN^{\RN}j \tau(t_{\ast}) \pi i }
\nonumber\\
&= 2 \pi r \delta_{jk} \vartheta_2 \Big( \cN^{\RN}j \tau(t_{\ast}); (\cN^{\RN})^2 \tau(t_{\ast}) \Big).
\end{align*}
This proves (\ref{eqn:ortho_BD1}) for $\RN=\CN$ and $\BCN$
with (\ref{eqn:ortho_const1}). \\
(iv) We prove (\ref{eqn:ortho_BD1}) for $\RN=\CNv$.
By (\ref{eqn:complex_conjugates}), LHS of (\ref{eqn:ortho_BD1})
is given by
\[
I^{\CNv}_{jk} = - \int_0^{\pi r} M^{\CNv}_j(x, t_{\ast}-t) M^{\CNv}_k(x, t) dx
= - I^{\CNv}_{jk, +} + I^{\CNv}_{jk, -},
\]
where
\begin{align*}
I^{\CNv}_{jk, +} 
&= \int_{-\pi r}^{\pi r} dx \, e^{2 \pi i (j+k-1) \xi(x)}
\vartheta_2 \Big( \cN^{\CNv}\{(j-1/2) \tau(t_{\ast}-t) + \xi(x) \};
(\cN^{\CNv})^2 \tau(t_{\ast}-t) \Big)
\nonumber\\
& \qquad \qquad \qquad \times
\vartheta_2 \Big( \cN^{\CNv}\{(k-1/2) \tau(t) + \xi(x) \};
(\cN^{\CNv})^2 \tau(t) \Big),
\nonumber\\
I^{\CNv}_{jk, -} 
&= \int_{-\pi r}^{\pi r} dx \, e^{2 \pi i (j-k) \xi(x)}
\vartheta_2 \Big( \cN^{\CNv}\{(j-1/2) \tau(t_{\ast}-t) + \xi(x) \};
(\cN^{\CNv})^2 \tau(t_{\ast}-t) \Big)
\nonumber\\
& \qquad \qquad \qquad \times
\vartheta_2 \Big( \cN^{\CNv}\{(k-1/2) \tau(t) - \xi(x) \};
(\cN^{\CNv})^2 \tau(t) \Big).
\end{align*}
We follow the similar argument to the case (iii).
Here the key inequalities are
$1 \leq j+k-1 \leq 2N-1 < \cN^{\CNv}=2N$
and
$0 \leq j-k \leq N-1 < \cN^{\CNv}$, $j, k \in \{1,2, \dots, N\}$.
Then we can conclude
$I^{\CNv}_{jk,+} \equiv 0$ and
\[
I^{\CNv}_{jk}
= I^{\CNv}_{jk, -}
= 2 \pi r \delta_{jk} \vartheta_2 \Big( \cN^{\CNv} (j-1/2) 
\tau(t_{\ast}); (\cN^{\CNv})^2 \tau(t_{\ast}) \Big).
\]
This proves (\ref{eqn:ortho_BD1}) for $\RN=\CNv$
with (\ref{eqn:ortho_const1}). \\
(v) Finally we prove (\ref{eqn:ortho_BD1}) for $\RN=\DN$.
By (\ref{eqn:complex_conjugates}), LHS of (\ref{eqn:ortho_BD1})
is given by
\[
I^{\DN}_{jk} = \int_0^{\pi r} M^{\DN}_j(x, t_{\ast}-t) M^{\DN}_k(x, t) dx
= I^{\DN}_{jk, +} + I^{\DN}_{jk, -},
\]
where
\begin{align*}
I^{\DN}_{jk, +} 
&= \int_{-\pi r}^{\pi r} dx \, e^{2 \pi i (j+k-2) \xi(x)}
\vartheta_2 \Big( \cN^{\DN}\{(j-1) \tau(t_{\ast}-t) + \xi(x) \};
(\cN^{\DN})^2 \tau(t_{\ast}-t) \Big)
\nonumber\\
& \qquad \qquad \qquad \times
\vartheta_2 \Big( \cN^{\DN}\{(k-1) \tau(t) + \xi(x) \};
(\cN^{\DN})^2 \tau(t) \Big),
\nonumber\\
I^{\DN}_{jk, -} 
&= \int_{-\pi r}^{\pi r} dx \, e^{2 \pi i (j-k) \xi(x)}
\vartheta_2 \Big( \cN^{\DN}\{(j-1) \tau(t_{\ast}-t) + \xi(x) \};
(\cN^{\DN})^2 \tau(t_{\ast}-t) \Big)
\nonumber\\
& \qquad \qquad \qquad \times
\vartheta_2 \Big( \cN^{\DN}\{(k-1) \tau(t) - \xi(x) \};
(\cN^{\DN})^2 \tau(t) \Big).
\end{align*}
By the definition of $\vartheta_2$ given by (\ref{eqn:theta}), we have
\begin{align*}
I^{\DN}_{jk, +}
&= \sum_{n \in \Z} \sum_{m \in \Z}
e^{ (n-1/2)^2 (\cN^{\DN})^2 \tau(t_{\ast}-t) \pi i
+(2n-1) \cN^{\DN}(j-1) \tau(t_{\ast}-t) \pi i }
\nonumber\\
& \qquad \times
e^{
(m-1/2)^2 (\cN^{\DN})^2 \tau(t) \pi i + (2m-1) \cN^{\DN} (k-1) \tau(t) \pi i }
\int_{-\pi r}^{\pi r} e^{2 \pi i \{(j+k-2)+\cN^{\DN}(n+m-1)\}\xi(x)} dx,
\nonumber\\
I^{\DN}_{jk, -}
&= \sum_{n \in \Z} \sum_{m \in \Z}
e^{ (n-1/2)^2 (\cN^{\DN})^2 \tau(t_{\ast}-t) \pi i
+(2n-1) \cN^{\DN}(j-1) \tau(t_{\ast}-t) \pi i }
\nonumber\\
& \qquad \times
e^{
(m-1/2)^2 (\cN^{\DN})^2 \tau(t) \pi i + (2m-1) \cN^{\DN} (k-1) \tau(t) \pi i }
\int_{-\pi r}^{\pi r} e^{2 \pi i \{(j-k)+\cN^{\DN}(n-m)\}\xi(x)} dx.
\end{align*}
Here we note that
\begin{align*}
\int_{-\pi r}^{\pi r} e^{2 \pi i \{(j+k-2)+\cN^{\DN}(n+m-1)\} \xi(x)} dx
&= 2 \pi r {\bf 1}\Big( (j+k-2)+\cN^{\DN}(n+m-1)=0 \Big),
\nonumber\\
\int_{-\pi r}^{\pi r} e^{2 \pi i \{(j-k)+\cN^{\DN}(n-m)\} \xi(x)} dx
&= 2 \pi r {\bf 1}\Big( (j-k)+\cN^{\DN}(n-m)=0 \Big).
\end{align*}
Since $j, k \in \{1,2, \dots, N\}$, we see that
\[
0 \leq j+k-2 \leq 2(N-1) = \cN^{\DN}, 
\]
and
$0 \leq j-k \leq N-1 < \cN^{\DN}$.
The condition $(j+k-2)+\cN^{\DN}(n+m-1)=0$ is satisfied, if
\[
j+k-2=0, \quad n+m-1=0
\quad \Longleftrightarrow \quad
j=k=1, \quad m=-n+1,
\]
or if
\[
j+k-2=2(N-1), \quad n+m-1=-1
\quad \Longleftrightarrow \quad
j=k=N, \quad m=-n.
\]
And the condition $(j-k)+\cN^{\DN}(n-m)=0$ is
satisfied if and only if
$j=k, n=m$.
Hence we see that
\[
I^{\DN}_{jk, +}=\delta_{j1} \delta_{k1} J_1
+\delta_{jN} \delta_{kN} J_2,
\]
where
$J_1 = 2 \pi r \vartheta_2 (0; (\cN^{\DN})^2 \tau(t_{\ast}) )$,
and
\begin{align*}
J_2& =
2 \pi r \sum_{n \in \Z} e^{
(n-1/2)^2 (\cN^{\DN})^2 \tau(t_{\ast}-t) \pi i 
+(2n-1) \cN^{\DN}(N-1) \tau(t_{\ast}-t) \pi i }
\nonumber\\
& \qquad \qquad
\times
e^{
(-n-1/2)^2 (\cN^{\DN})^2 \tau(t) \pi i 
+(-2n-1) \cN^{\DN}(N-1) \tau(t) \pi i }.
\end{align*}
By the fact $\cN^{\DN}=2(N-1)$, 
it is easy to verify that 
\begin{align*}
&
(n-1/2)^2 (\cN^{\DN})^2 \tau(t_{\ast}-t)
+(2n-1) \cN^{\DN}(N-1) \tau(t_{\ast}-t) 
\nonumber\\
& \quad +(-n-1/2)^2 (\cN^{\DN})^2 \tau(t) 
+(-2n-1) \cN^{\DN}(N-1) \tau(t) 
\nonumber\\
& \qquad =(n-1/2)^2 (\cN^{\DN})^2 \tau(t_{\ast})
+(2n-1) \cN^{\DN}(N-1) \tau(t_{\ast}),
\end{align*}
and thus
$J_2
=2 \pi r \vartheta_2 ( \cN^{\DN}(N-1) \tau(t_{\ast}); 
(\cN^{\DN})^2 \tau(t_{\ast}) )$.
We also obtain 
\[
I^{\DN}_{jk,-}
= 2 \pi r \delta_{jk} \vartheta_2 ( \cN^{\DN}(j-1) \tau(t_{\ast}); 
(\cN^{\DN})^2 \tau(t_{\ast}) ).
\]
Therefore, we can conclude
\begin{align*}
I^{\DN}_{jk} &= I^{\DN}_{jk, +}+ I^{\DN}_{jk, -}
\nonumber\\
&= 2 \pi r \delta_{jk}
\Big\{ \delta_{j1}
\vartheta_2 \Big( 0; (\cN^{\DN})^2 \tau(t_{\ast}) \Big)
+
\vartheta_2 \Big( \cN^{\DN}(j-1) \tau(t_{\ast}); (\cN^{\DN})^2 \tau(t_{\ast}) \Big)
\nonumber\\
& \qquad \qquad + \delta_{jN}
\vartheta_2 \Big( \cN^{\DN}(N-1) \tau(t_{\ast}); (\cN^{\DN})^2 \tau(t_{\ast}) \Big)
\Big\}.
\end{align*}
This proves (\ref{eqn:ortho_BD1}) for $\RN=\DN$
with (\ref{eqn:ortho_const3}).
The proof is complete. \qed

\vskip 0.3cm
\noindent{\bf Remark 1} \,
When $t=t_{\ast}/2$, the functions
$\{M^{\RN}_j(x,t_{\ast}/2)\}_{j=1}^N$ form orthogonal sets
with respect to the inner product
$\langle f | g \rangle = \int_0^{L} \overline{f(x)} g(x) dx$
with $L=2 \pi r$ for $\RN=\AN$ 
and $L=\pi r$ for $\RN=\BN$,$\BNv$,$\CN$,$\CNv$,$\BCN$,$\DN$.
For the case $\RN=\AN$, this fact was announced
on page 217 in \cite{For10}.
\vskip 0.3cm

\SSC
{Determinantal Point Processes} \label{sec:DPP}
\subsection{Main results} \label{sec:main_results}

As functions of $\tau$, we define
\[
q(\tau) = e^{\tau \pi i},
\qquad 
q_0(\tau) =
\prod_{n=1}^{\infty} (1-q(\tau)^{2n}). 
\]
In our setting (\ref{eqn:MR0}) with (\ref{eqn:xi_tau1}), 
the Macdonald denominator formulas
(\ref{eqn:Macdonald1}) and (\ref{eqn:Macdonald2}) of
Rosengren and Schlosser (Proposition 6.1 in \cite{RS06})
are written as follows.
(For type $\AN$, we set the norm as 
$\alpha=e^{2 \pi i \widetilde{\alpha}_N}$ with
(\ref{eqn:norm}), which is different from
the choice in Lemma 2.4 in \cite{Kat17}).
\begin{align}
& \det_{1 \leq j, k \leq N}
\Big[ M^{\AN}_j(x_k, t) \Big]
\nonumber\\
& \quad
=\begin{cases}
\displaystyle{
i^{N/2} a^{\AN}(t)
\vartheta_0 \left(\sum_{j=1}^N \xi(x_{j}) ;
\cN^{\AN} \tau(t) \right)
W^{\AN}(\xi(\x); \cN^{\AN} \tau(t))
},
& \mbox{if $N$ is even},
\cr
\displaystyle{
i^{-(N-1)/2} a^{\AN}(t)
\vartheta_3 \left(\sum_{j=1}^N \xi(x_{j}) ;
\cN^{\AN} \tau(t) \right)
W^{\AN}(\xi(\x); \cN^{\AN} \tau(t))
},
& \mbox{if $N$ is odd},
\end{cases}
\nonumber\\
& \det_{1 \leq j, k \leq N}
\Big[ M^{\RN}_j(x_k, t) \Big]
= a^{\RN}(t)
W^{\RN}(\xi(\x); \cN^{\RN} \tau(t)),
\quad \mbox{for $\RN=\BN, \BNv, \DN$},
\nonumber\\
& \det_{1 \leq j, k \leq N}
\Big[ M^{\RN}_j(x_k, t) \Big]
= i^{-N} a^{\RN}(t)
W^{\RN}(\xi(\x); \cN^{\RN} \tau(t)),
\quad \mbox{for $\RN=\CN, \CNv, \BCN$},
\label{eqn:Macdonald3}
\end{align}
where
\begin{align}
a^{\AN}(t) &= q(\cN^{\AN} \tau(t))^{-N(3N-1)/8} 
q_0(\cN^{\AN} \tau(t))^{-(N-1)(N-2)/2},
\nonumber\\
a^{\BN}(t) &= 2 q(\cN^{\BN} \tau(t))^{-N(N-1)/4} 
q_0(\cN^{BN} \tau(t))^{-N(N-1)},
\nonumber\\
a^{\BNv}(t) &= 2 q(\cN^{\BNv} \tau(t))^{-N(N-1)/4} 
q_0(\cN^{\BNv} \tau(t) )^{-(N-1)^2} q_0(2 \cN^{\BNv} \tau(t))^{-(N-1)},
\nonumber\\
a^{\CN}(t) &= q(\cN^{\CN} \tau(t))^{-N^2/4} 
q_0(\cN^{\CN} \tau(t))^{-N(N-1)},
\nonumber\\
a^{\CNv}(t) &= q(\cN^{\CNv} \tau(t))^{-N(2N-1)/8} 
q_0(\cN^{\CNv} \tau(t))^{-(N-1)^2} q_0(\cN^{\CNv} \tau(t)/2)^{-(N-1)},
\nonumber\\
a^{\BCN}(t) &= q(\cN^{\BCN} \tau(t))^{-N(N+1)/4} 
q_0(\cN^{\BCN} \tau(t) )^{-N(N-1)} q_0(2 \cN^{\BCN} \tau(t))^{-N},
\nonumber\\
a^{\DN}(t) &= 4 q(\cN^{\DN} \tau(t) )^{-N(N-1)/4} 
q_0(\cN^{\DN} \tau(t))^{-N(N-2)},
\label{eqn:a1}
\end{align}
and
$\xi(\x) \equiv (\xi(x_1), \xi(x_2), \dots, \xi(x_N))$.
Note that
\begin{align*}
& q(\cN^{\RN} \tau(t)) = e^{-\cN^{\RN} t/2r^2} > 0,
\nonumber\\
& q_0(\cN^{\RN} \tau(t)) = 
\prod_{n=1}^{\infty}(1-e^{-n \cN^{\RN} t/r^2}) \geq 0,
\quad \mbox{if $0 \leq t < \infty$}.
\end{align*}

Consider the following Weyl alcoves,
\begin{align*}
\W_N^{[0, 2 \pi r)}
&\equiv \{\x =(x_1, x_2, \dots, x_N) \in \R^N: 
0 \leq x_1< x_2 < \cdots < x_N < 2 \pi r \},
\nonumber\\
\W_N^{[0, \pi r]}
&\equiv \{\x =(x_1, x_2, \dots, x_N) \in \R^N: 
0 \leq x_1< x_2 < \cdots < x_N \leq \pi r \}.
\end{align*}
By (\ref{eqn:theta_positive}),
$
\vartheta_s ( \sum_{j=1}^N \xi(x_j);
\cN^{\AN} \tau(t)) \geq 0$
for $s=0,3$, if $t \geq 0$, 
and the definitions of Macdonald denominators
(\ref{eqn:Macdonald_denominators}) imply that
\begin{align*}
& W^{\AN}(\xi(\x); \cN^{\AN} \tau(t)) \geq 0,
\quad 
\mbox{if $\x \in \W_N^{[0, 2 \pi r)}, \, t \geq 0$},
\nonumber\\
& W^{\RN}(\xi(\x); \cN^{\RN} \tau(t)) \geq 0,
\quad 
\mbox{if $\x \in \W_N^{[0, \pi r]}, \, t \geq 0$},
\, \mbox{for $\RN=\BN,\BNv,\CN,\CNv,\BCN,\DN$}.
\end{align*}
Now we introduce
\begin{equation}
q^{\RN}_t(\x)
=\det_{1 \leq j , k \leq N}
\left[ \overline{M^{\RN}_{j}(x_k, t_{\ast}-t)} \right]
\det_{1 \leq \ell, m \leq N}
\left[ M^{\RN}_{\ell}(x_m, t) \right],
\quad t \in (0, t_{\ast}). 
\label{eqn:q1}
\end{equation}
By the basic properties of
the Jacobi theta functions (\ref{eqn:even_odd})--(\ref{eqn:theta_positive}),
the product form of (\ref{eqn:q1}) guarantees the following.

\begin{lem}
\label{thm:positivity}
If $t \in (0, t_{\ast})$, 
$q^{\RN}_t(\x) \geq 0$, $\x \in \R^N$, for 
$\RN=\AN,\BN,\BNv,\CN,\CNv,\BCN,\DN$.
\end{lem}
\vskip 0.3cm
Moreover, we can verify the following.
\begin{lem}
\label{thm:normalization}
For $t \in (0, t_{\ast})$, 
\begin{align}
\int_{\W_N^{[0, 2 \pi r)}}
q^{\AN}_t(\x) d \x
&= \prod_{n=1}^N m^{\AN}_n(t_{\ast}),
\label{eqn:normalization1}
\\
\int_{\W_N^{[0, \pi r]}}
q^{\RN}_t(\x) d \x
&= \prod_{n=1}^N m^{\RN}_n(t_{\ast}),
\quad 
\RN=\BN,\BNv,\CN,\CNv,\BCN,\DN.
\label{eqn:normalization2}
\end{align}
\end{lem}
\noindent{\it Proof} \,
Let $S^{\AN}=\W_N^{[0, 2 \pi r)}$,
$L^{\AN}=2 \pi r$, and
$S^{\RN}=\W_N^{[0, \pi r]}$,
$L^{\RN}= \pi r$ for 
$\RN=\BN$, $\BNv$, $\CN$, $\CNv$, $\BCN$, $\DN$.
By the Heine identity (\ref{eqn:Heine}) in Appendix \ref{sec:appendixC}, 
which is also called the Andr\'eief or Gram identity, 
\[
\int_{S^{\RN}} q^{\RN}_t(\x) d \x
=\det_{1 \leq j, k \leq N}
\left[ \int_0^{L^{\RN}} \overline{M^{\RN}_j(x, t_{\ast}-t)}
M^{\RN}_k(x,t) dx \right].
\]
By the biorthogonality given by Lemma \ref{thm:orthogonality},
this is equal to 
$\det_{1 \leq j,k \leq N} [ m^{\RN}_j(t_{\ast}) \delta_{jk} ]$,
and hence (\ref{eqn:normalization1}) and (\ref{eqn:normalization2})
are proved. \qed
\vskip 0.3cm
\vskip 0.3cm
\noindent{\bf Remark 2} \,
Combining this lemma with the Macdonald denominator formulas
(\ref{eqn:Macdonald3}) with (\ref{eqn:a1}), 
we readily obtain the Selberg-type integral formulas 
(see, for instance, Chapter 14 of \cite{For10})
for products of Macdonald denominators, 
see Appendix \ref{sec:appendixB}.
They seem to be much simpler than the
formulas known as
{\it elliptic Selberg integrals} 
(see Section 4.4 of \cite{FW08} , 
Exercise 4.1.4 in \cite{For10},
and references therein). 
\vskip 0.3cm

Then the seven types of
one-parameter ($t \in (0, t_{\ast})$) families 
of probability measures $\bP^{\RN}_t$ are defined as
\begin{align}
\bP^{\RN}_t(\X \in d \x)
=\bp^{\RN}_t(\x) d \x
&= \begin{cases}
\displaystyle{
\frac{q^{\AN}_t(\x)}{\prod_{n=1}^N m_n(t_{\ast})} d \x
}, & 
\mbox{for $\RN=\AN$}, 
\cr
& \cr
\displaystyle{
\frac{q^{\RN}_t(\x)}{\prod_{n=1}^N m_n(t_{\ast})} d \x
}, &
\mbox{for $\RN=\BN,\BNv,\CN,\CNv,\BCN,\DN$},
\end{cases}
\label{eqn:PtR1}
\end{align}
which are normalized as
\begin{align}
& \int_{\W^{[0, 2 \pi r)}} \bp^{\AN}(\x) d \x=1,
\nonumber\\
& \int_{\W^{[0, \pi r]}} \bp^{\RN}(\x) d \x =1,
\quad \RN=\BN, \BNv, \CN, \CNv, \BCN, \DN.
\label{eqn:normalization_p}
\end{align}
Under these probability measures 
$\bP^{\RN}_t$ with one parameter $t \in (0, t_{\ast})$, 
we consider seven types of point processes, 
\[
\Xi^{\AN}(\cdot)=\sum_{j=1}^N \delta_{X^{\AN}_j}(\cdot) \quad
\mbox{on $S=[0, 2 \pi r)$}, 
\]
and
\[
\Xi^{\RN}(\cdot)=\sum_{j=1}^N \delta_{X^{\RN}_j}(\cdot) \quad
\mbox{on $S=[0, \pi r]$}, \quad
\mbox{for $\RN=\BN, \BNv, \CN, \CNv, \BCN, \DN$}.
\]
Given the determinantal expressions
(\ref{eqn:q1}) and (\ref{eqn:PtR1}) for the
probability measures associated with
the biorthogonal functions
(\ref{eqn:MA1})--(\ref{eqn:MD1}), we can readily prove the following fact
by the standard method 
in random matrix theory \cite{Meh04,For10,AGZ10,Kat15_Springer}. 
We give a sketch of proof for a general statement 
in Appendix \ref{sec:appendixC}
for convenience of the reader.

\begin{thm}
\label{thm:mainA1}
The seven types of 
one-parameter families of point processes, 
$(\Xi^{\RN}, \bP_t^{\RN}, t \in (0, t_{\ast}))$,
$\RN=\AN$, $\BN$, $\BNv$, $\CN$, $\CNv$, $\BCN$, $\DN$, 
are determinantal with the correlation kernels, 
\begin{align}
K_t^{\RN}(x,y; t_{\ast}, r)  &= \sum_{n=1}^{N}
\frac{1}{m^{\RN}_n(t_{\ast})} M^{\RN}_n(x, t) 
\overline{M^{\RN}_n(y, t_{\ast}-t)}, 
\quad t \in (0, t_{\ast}), 
\nonumber\\
&
x, y \in [0, 2\pi r), \quad \mbox{for $\RN=\AN$},
\nonumber\\
&
x, y \in [0, \pi r], \quad 
\mbox{for $\RN=\BN,\BNv,\CN,\CNv,\BCN,\DN$}.
\label{eqn:K_R1}
\end{align}
\end{thm}

\subsection{Temporally homogeneous limit at $t=t_{\ast}/2$}
\label{sec:temp_homo_finite}

We consider the determinantal point processes at $t=t_{\ast}/2$.
The correlation kernels (\ref{eqn:K_R1}) become
\begin{equation}
K_{t_{\ast}/2}^{\RN}(x,y; t_{\ast}, r)  = \sum_{n=1}^{N}
\frac{1}{m^{\RN}_n(t_{\ast})} M^{\RN}_n(x, t_{\ast}/2) 
\overline{M^{\RN}_n(y, t_{\ast}/2)}, 
\label{eqn:K_R2}
\end{equation}
$x, y \in [0, 2\pi r)$ for $\RN=\AN$,
and $x, y \in [0, \pi r]$ 
for $\RN=\BN$, $\BNv$, $\CN$, $\CNv$, $\BCN$, $\DN$.

By the asymptotics of the Jacobi theta functions (\ref{eqn:theta_asym}),
the temporally homogeneous limit $t_{\ast} \to \infty$ of
(\ref{eqn:K_R2}) are obtained as follows. \\
\noindent
(i) For $\RN=\AN$, 
\begin{align}
K^{\AN}(x,y; r) 
&\equiv \lim_{t_{\ast} \to \infty} K^{\AN}_{t_{\ast}/2}(x,y;t_{\ast}, r)
\nonumber\\
&= \frac{1}{2 \pi r}
\sum_{n=1}^N e^{2 \pi i (n-1) (\xi(x)-\xi(y))}
= 
\frac{1}{2 \pi r}
\frac{\sin \{N(x-y)/2r\}}{\sin \{(x-y)/2r\}},
\quad x, y \in [0, 2 \pi r).
\label{eqn:K_f_h_A2}
\end{align}
(ii) For $\RN=\BN, \BNv, \CN, \CNv, \BCN$,
\begin{align*}
& K^{\RN}(x,y;r)
\equiv
\lim_{t_{\ast} \to \infty} K^{\RN}_{t_{\ast}/2}(x,y; t_{\ast}, r)
\nonumber\\
& \quad
= \frac{2}{\pi r}
\sum_{n=1}^N \sin \{\pi(\cN^{\RN}-2 J^{\RN}(n)) \xi(x) \}
\sin \{\pi(\cN^{\RN}-2 J^{\RN}(n)) \xi(y) \}
\nonumber\\
& \quad
= \begin{cases}
\displaystyle{\frac{1}{2 \pi r} \left[
\frac{\sin\{(\cN^{\RN}+1)(x-y)/2r\}}{\sin\{(x-y)/2r\}}
- \frac{\sin\{(\cN^{\RN}+1)(x+y)/2r\}}{\sin\{(x+y)/2r\}}
\right]},
& \mbox{if $\RN=\BN, \BNv$}, 
\cr
& \cr
\displaystyle{\frac{1}{2 \pi r} \left[
\frac{\sin\{(\cN^{\RN}-1)(x-y)/2r\}}{\sin\{(x-y)/2r\}}
- \frac{\sin\{(\cN^{\RN}-1)(x+y)/2r\}}{\sin\{(x+y)/2r\}}
\right]},
& \mbox{if $\RN=\CN, \BCN$}, 
\cr
& \cr
\displaystyle{\frac{1}{2 \pi r} \left[
\frac{\sin\{\cN^{\RN}(x-y)/2r\}}{\sin\{(x-y)/2r\}}
- \frac{\sin\{\cN^{\RN}(x+y)/2r\}}{\sin\{(x+y)/2r\}}
\right]},
& \mbox{if $\RN=\CNv$},
\cr
\end{cases}
\end{align*}
$x, y \in [0, \pi r]$. \\
(iii) For $\RN=\DN$,
\begin{align}
& K^{\DN}(x,y;r)
\equiv
\lim_{t_{\ast} \to \infty} K^{\DN}_{t_{\ast}/2}(x,y; t_{\ast}, r)
\nonumber\\
& \quad
= \frac{2}{\pi r}
\sum_{n=1}^N \cos \{2\pi(N-n) \xi(x) \}
\cos \{2\pi(N-n) \xi(y) \}
\nonumber\\
& \quad
= \frac{1}{2 \pi r} \left[
\frac{\sin\{(2N-1)(x-y)/2r\}}{\sin\{(x-y)/2r\}}
+ \frac{\sin\{(2N-1)(x+y)/2r\}}{\sin\{(x+y)/2r\}}
\right], \quad x, y \in [0, \pi r]. 
\label{eqn:K_f_h_D}
\end{align}
Since $\cN^{\BN}+1=\cN^{\BCN}-1=\cN^{\CNv}=2N$,
and
$\cN^{\BNv}+1=\cN^{\CN}-1=2N+1$, 
\begin{align}
& K^{\BN}(x,y;r) =K^{\BCN}(x,y;r)=K^{\CNv}(x,y;r)
\nonumber\\
& \quad =
\frac{1}{2 \pi r} \left[
\frac{\sin\{N(x-y)/r\}}{\sin\{(x-y)/2r\}}
- \frac{\sin\{N(x+y)/r\}}{\sin\{(x+y)/2r\}}
\right], \quad x, y \in [0, \pi r], 
\label{eqn:K_f_h_B}
\\
& K^{\CN}(x,y;r) = K^{\BNv}(x,y;r) 
\nonumber\\
& \quad =
\frac{1}{2 \pi r} \left[
\frac{\sin\{(2N+1)(x-y)/2r\}}{\sin\{(x-y)/2r\}}
- \frac{\sin\{(2N+1)(x+y)/2r\}}{\sin\{(x+y)/2r\}}
\right], \quad x, y \in [0, \pi r].
\label{eqn:K_f_h_C}
\end{align}

\begin{cor}
\label{thm:finite_homo}
Put $t=t_{\ast}/2$ in Theorem \ref{thm:mainA1}. 
In the limit $t_{\ast} \to \infty$, 
the seven types of determinantal point processes
$(\Xi^{\RN}, \bP^{\RN}_{t_{\ast}/2})$
are degenerated into the four types of determinantal point processes 
specified by the correlation kernels 
$K^{\AN}(x, y; r)$, $K^{\BN}(x,y;r)$,
$K^{\CN}(x,y;r)$, and $K^{\DN}(x,y;r)$
as shown by 
{\rm (\ref{eqn:K_f_h_A2})}, {\rm (\ref{eqn:K_f_h_B})}, 
{\em (\ref{eqn:K_f_h_C})}, and {\rm (\ref{eqn:K_f_h_D})},
respectively.
\end{cor}

\vskip 0.3cm
\noindent{\bf Remark 3} \,
The correlation kernel $K^{\AN}(x, y; r)$ determines the equilibrium determinantal
point processes of the noncolliding BMs
on a circle with radius $r>0$ 
(see \cite{NF03} and Proposition 6.1 in \cite{Kat14}),
and the correlation kernels $K^{\CN}(x, y; r)$ and $K^{\DN}(x,y; r)$ 
do the equilibrium point processes 
in an interval $[0, \pi r]$ with the
absorbing and reflecting boundary conditions, respectively
(see Proposition 5.2 in \cite{Kat17}). 
In random matrix theory, the determinantal point processes
governed by the correlation kernels
$K^{\RN}(x, y; 1)$ of the four types,
$\RN=\AN, \BN, \CN, \DN$ are realized as the 
eigenvalue distributions of random matrices
in ${\rm U}(N)$, ${\rm SO}(2N+1)$,
${\rm Sp}(N)$, ${\rm SO(2N)}$, respectively, 
see Section 2.3 (c) in \cite{Sos00}.
In particular, the eigenvalue distribution
of random matrices in ${\rm U}(N)$ is called the
{\it circular unitary ensemble} (CUE), 
see Chapter 11 in \cite{Meh04}.
\vskip 0.3cm
\subsection{Infinite determinantal point processes}
\label{sec:infinite}

We fix the density of points as
\begin{equation}
\rho=\begin{cases}
\displaystyle{\frac{N}{2 \pi r}},
& \RN=\AN,
\cr
& \cr
\displaystyle{ \frac{N}{\pi r}},
&\RN=\BN, \BNv, \CN, \CNv, \BCN, \DN,
\end{cases}
\label{eqn:rho1}
\end{equation}
and take double limit $N \to \infty$, $r \to \infty$.
Then we obtain the following limits of correlation kernels.

\begin{lem}
\label{thm:scaling_limit_K}
For $t \in (0, t_{\ast})$, 
the following scaling limits are obtained for correlation kernels. \\
{\rm (i)} For $\RN=\AN$, 
\begin{align}
& \cK^{A}_t(x, y; t_{\ast}, \rho)
\equiv \lim_{\substack{N \to \infty, r \to \infty, \cr N/2\pi r=\rho}}
K^{\AN}_t(x, y; t_{\ast}, r)
\nonumber\\
& \quad
= \int_0^{\rho} d \lambda \,
e^{2\pi i (x-y) \lambda} 
\frac{\vartheta_2 (\rho x+2 \pi i t \rho \lambda; 2 \pi i t \rho^2)
\vartheta_2 (\rho y - 2 \pi i (t_{\ast}-t) \rho \lambda; 2 \pi i (t_{\ast}-t) \rho^2)}
{\vartheta_2 (2 \pi i t_{\ast} \rho \lambda; 2  \pi i t_{\ast} \rho^2)},
\label{eqn:K_inf_A}
\end{align}
$x, y \in \R$. \\
{\rm (ii)} For $\RN=\BN, \BNv$,
\begin{align}
& \cK^{B}_t(x, y; t_{\ast}, \rho)
\equiv \lim_{\substack{N \to \infty, r \to \infty, \cr N/\pi r=\rho}}
K^{\RN}_t(x, y; t_{\ast}, r)
\nonumber\\
& \quad
= \frac{1}{2} \left[
\int_{-\rho}^{\rho} d \lambda \,
e^{\pi i (x-y) \lambda} 
\frac{\vartheta_1 (\rho x+\pi i t \rho \lambda; 2 \pi i t \rho^2)
\vartheta_1 (\rho y - \pi i (t_{\ast}-t) \rho \lambda; 2 \pi i (t_{\ast}-t) \rho^2)}
{\vartheta_2 (\pi i t_{\ast} \rho \lambda; 2  \pi i t_{\ast} \rho^2)} 
\right.
\nonumber\\
& \qquad \left. \quad
- \int_{-\rho}^{\rho} d \lambda \,
e^{\pi i (x+y) \lambda} 
\frac{\vartheta_1 (\rho x+\pi i t \rho \lambda; 2 \pi i t \rho^2)
\vartheta_1 (-\rho y - \pi i (t_{\ast}-t) \rho \lambda; 2 \pi i (t_{\ast}-t) \rho^2)}
{\vartheta_2 (\pi i t_{\ast} \rho \lambda; 2  \pi i t_{\ast} \rho^2)} \right],
\label{eqn:K_inf_B}
\end{align}
$x, y \in [0, \infty)$. \\
{\rm (iii)} For $\RN=\CN, \CNv, \BCN$,
\begin{align}
& \cK^{C}_t(x, y; t_{\ast}, \rho)
\equiv \lim_{\substack{N \to \infty, r \to \infty, \cr N/\pi r=\rho}}
K^{\RN}_t(x, y; t_{\ast}, r)
\nonumber\\
& \quad
= \frac{1}{2} \left[
\int_{-\rho}^{\rho} d \lambda \,
e^{\pi i (x-y) \lambda} 
\frac{\vartheta_2 (\rho x+\pi i t \rho \lambda; 2 \pi i t \rho^2)
\vartheta_2 (\rho y - \pi i (t_{\ast}-t) \rho \lambda; 2 \pi i (t_{\ast}-t) \rho^2)}
{\vartheta_2 (\pi i t_{\ast} \rho \lambda; 2  \pi i t_{\ast} \rho^2)} 
\right.
\nonumber\\
& \qquad \left. \quad
- \int_{-\rho}^{\rho} d \lambda \,
e^{\pi i (x+y) \lambda} 
\frac{\vartheta_2 (\rho x+\pi i t \rho \lambda; 2 \pi i t \rho^2)
\vartheta_2 (-\rho y - \pi i (t_{\ast}-t) \rho \lambda; 2 \pi i (t_{\ast}-t) \rho^2)}
{\vartheta_2 (\pi i t_{\ast} \rho \lambda; 2  \pi i t_{\ast} \rho^2)} \right],
\label{eqn:K_inf_C}
\end{align}
$x, y \in [0, \infty)$. \\
{\rm (iv)} For $\RN=\DN$,
\begin{align}
& \cK^{D}_t(x, y; t_{\ast}, \rho)
= \lim_{\substack{N \to \infty, r \to \infty, \cr N/\pi r=\rho}}
K^{\DN}_t(x, y; t_{\ast}, r)
\nonumber\\
& \quad
= \frac{1}{2} \left[
\int_{-\rho}^{\rho} d \lambda \,
e^{\pi i (x-y) \lambda} 
\frac{\vartheta_2 (\rho x+\pi i t \rho \lambda; 2 \pi i t \rho^2)
\vartheta_2 (\rho y - \pi i (t_{\ast}-t) \rho \lambda; 2 \pi i (t_{\ast}-t) \rho^2)}
{\vartheta_2 (\pi i t_{\ast} \rho \lambda; 2  \pi i t_{\ast} \rho^2)} 
\right.
\nonumber\\
& \qquad \left. \quad
+ \int_{-\rho}^{\rho} d \lambda \,
e^{\pi i (x+y) \lambda} 
\frac{\vartheta_2 (\rho x+\pi i t \rho \lambda; 2 \pi i t \rho^2)
\vartheta_2 (-\rho y - \pi i (t_{\ast}-t) \rho \lambda; 2 \pi i (t_{\ast}-t) \rho^2)}
{\vartheta_2 (\pi i t_{\ast} \rho \lambda; 2  \pi i t_{\ast} \rho^2)} \right],
\label{eqn:K_inf_D}
\end{align}
$x, y \in [0, \infty)$. \\
\end{lem}
\noindent{\it Proof} \,
Here we give proof for (iii).
Other cases are similarly proved.
The explicit expressions for $K^{\RN}_t(x,y; t_{\ast}, r)$
for $\RN=\CN, \CNv, \BCN$ are given by
\begin{align}
&K^{\RN}_t(x, y; t_{\ast}, r) 
= - \frac{1}{2 \pi r} \sum_{n=1}^{N}
\frac{1}{\vartheta_2 (\cN^{\RN} J^{\RN}(n) \tau(t_{\ast}); 
(\cN^{\RN})^2 \tau(t_{\ast})) }
\nonumber\\
& \qquad \qquad
\times\Big\{ e^{2 \pi i J^{\RN}(n) \xi(x)}
\vartheta_2(\cN^{\RN}\{J^{\RN}(n) \tau(t)+\xi(x)\}; (\cN^R)^2 \tau(t))
\nonumber\\
& \qquad \qquad \qquad
- e^{-2 \pi i J^{\RN}(n) \xi(x)}
\vartheta_2(\cN^{\RN}\{J^{\RN}(n) \tau(t)-\xi(x)\}; (\cN^R)^2 \tau(t)) \Big\}
\nonumber\\
& \qquad \qquad
\times\Big\{ e^{2 \pi i J^{\RN}(n) \xi(y)}
\vartheta_2(\cN^{\RN}\{J^{\RN}(n) \tau(t_{\ast}-t)+\xi(y)\}; 
(\cN^R)^2 \tau(t_{\ast}-t))
\nonumber\\
& \qquad \qquad \qquad
- e^{-2 \pi i J^{\RN}(n) \xi(y)}
\vartheta_2(\cN^{\RN}\{J^{\RN}(n) \tau(t_{\ast}-t)-\xi(y)\}; 
(\cN^R)^2 \tau(t_{\ast}-t)) \Big\},
\label{eqn:K_limit1}
\end{align}
where $\cN^{\CN}=2(N+1)$, $\cN^{\CNv}=2N$, $\cN^{\BCN}=2N+1$,
and $J^{\CN}(n)=J^{\BCN}(n)=n$, $J^{\CNv}(n)=n-1/2$.
By (\ref{eqn:xi_tau1}) and (\ref{eqn:rho1}), we see that
\begin{align*}
& \frac{1}{2\pi r}=\frac{\rho}{2N}, \quad
\cN^{\RN} J^{\RN}(n) \tau(t)
=\frac{\cN^{\RN}}{2N} \pi i t \rho^2 \frac{J^{\RN}(n)}{N}, \quad
(\cN^{\RN})^2 \tau(t)= \left( \frac{\cN^{\RN}}{2N} \right)^2 2 \pi i t \rho^2,
\nonumber\\
& 
2 \pi i J^{\RN}(n) \xi(x)=\pi x \rho \frac{J^{\RN}(n)}{N},
\quad
\cN^{\RN} \xi(x)=\frac{\cN^{\RN}}{2N} \rho x.
\end{align*}
Since $N^{\RN}/2N \to 1$ as $N \to \infty$, 
(\ref{eqn:K_limit1}) given by summation converges 
uniformly on any compact subset of $[0, \pi r]^2 \ni (x,y)$ 
to the 
following integral with an integral variable $u \sim J^{\RN}(n)/N$,
\begin{align}
& - \frac{\rho}{2} \int_0^1 du \,
\frac{1}{\vartheta_2(\pi i t_{\ast} \rho^2 u; 2 \pi i t_{\ast} \rho^2)}
\Big\{ e^{\pi i x \rho u} \vartheta_2(\pi i t \rho^2 u+ \rho x; 2 \pi i t \rho^2)
- e^{-\pi i x \rho u} \vartheta_2(\pi i t \rho^2 u- \rho x; 2 \pi i t \rho^2) \Big\}
\nonumber\\
& \qquad \times
\Big\{ e^{\pi i y \rho u} \vartheta_2(\pi i (t_{\ast}-t) \rho^2 u+ \rho y; 
2 \pi i (t_{\ast}-t) \rho^2)
- e^{-\pi i y \rho u} \vartheta_2(\pi i (t_{\ast}-t) \rho^2 u- \rho y; 
2 \pi i (t_{\ast}-t) \rho^2) \Big\}.
\label{eqn:K_limit2}
\end{align}
We change the integral variable $u \to \lambda$ by $\lambda=\rho u$.
If we use the symmetry of Jacobi's theta functions, (\ref{eqn:even_odd}), 
we can verify that (\ref{eqn:K_limit2}) is rewritten as
(\ref{eqn:K_inf_C}). \qed
\vskip 0.3cm

The uniform convergence of correlation kernels implies 
the convergence of all correlation functions.
Then we conclude the following.

\begin{thm}
\label{thm:infinite_systems}
In the scaling limit $N \to \infty$, $r \to \infty$
with constant density of points {\rm (\ref{eqn:rho1})},
the seven types of 
one-parameter families of determinantal point processes, 
$(\Xi^{\RN}, \bP^{\RN}_t, t \in (0, t_{\ast}))$,
$\RN=\AN$, $\BN$, $\BNv$, $\CN$, $\CNv$, $\BCN$, $\DN$,
converge in the sense of finite dimensional distributions to
the four types of infinite dimensional point processes
as follows, 
\begin{align*}
(\Xi^{\AN}, \bP^{\AN}_t, t \in (0, t_{\ast}))
&\Longrightarrow
(\Xi^{A}, \bP^{A}_t, t \in (0, t_{\ast}))
\quad \mbox{
as $N \to \infty, r \to \infty$ with $\displaystyle{\frac{N}{2 \pi r} = \rho}$
},
\nonumber\\
\left.
\begin{array}{l}
(\Xi^{\BN}, \bP^{\BN}_t, t \in (0, t_{\ast})) \cr
(\Xi^{\BNv}, \bP^{\BNv}_t, t \in (0, t_{\ast})) 
\end{array}
\right\}
&\Longrightarrow
(\Xi^{B}, \bP^{B}_t, t \in (0, t_{\ast}))
\quad \mbox{
as $N \to \infty, r \to \infty$ with $\displaystyle{\frac{N}{\pi r} = \rho}$
},
\nonumber\\
\left.
\begin{array}{l}
(\Xi^{\CN}, \bP^{\CN}_t, t \in (0, t_{\ast})) \cr
(\Xi^{\CNv}, \bP^{\CNv}_t, t \in (0, t_{\ast})) \cr
(\Xi^{\BCN}, \bP^{\BCN}_t, t \in (0, t_{\ast})) 
\end{array}
\right\}
&\Longrightarrow
(\Xi^{C}, \bP^{C}_t, t \in (0, t_{\ast}))
\quad \mbox{as} \quad N \to \infty, r \to \infty, \frac{N}{\pi r} = \rho,
\nonumber\\
(\Xi^{\DN}, \bP^{\DN}_t, t \in (0, t_{\ast}))
&\Longrightarrow
(\Xi^{D}, \bP^{D}_t, t \in (0, t_{\ast}))
\quad \mbox{
as $N \to \infty, r \to \infty$ with $\displaystyle{\frac{N}{\pi r} = \rho}$
},
\end{align*}
where $(\Xi^{A}, \bP^{A}_t, t \in (0, t_{\ast}))$,
$(\Xi^{B}, \bP^{B}_t, t \in (0, t_{\ast}))$, 
$(\Xi^{C}, \bP^{C}_t, t \in (0, t_{\ast}))$, and
$(\Xi^{D}, \bP^{D}_t, t \in (0, t_{\ast}))$
are infinite determinantal point processes
associated with the correlation kernels
$\cK^{A}_t$, 
$\cK^{B}_t$, $\cK^{C}_t$, and $\cK^{D}_t$, $t \in (0, t_{\ast})$, 
which are given by {\rm (\ref{eqn:K_inf_A})},
{\rm (\ref{eqn:K_inf_B})},
{\rm (\ref{eqn:K_inf_C})}, and
{\rm (\ref{eqn:K_inf_D})},
respectively.
\end{thm}
\vskip 0.3cm
\noindent{\bf Remark 4} \,
Using (\ref{eqn:Theta}), we define
\begin{align*}
g^{A}(x, \lambda; t) &=
\frac{\Theta^{A}(\lambda/\rho, \rho x, 2 \pi i t \rho^2)}
{\sqrt{\Theta^{A}(\lambda/\rho, 0, 2 \pi i t_{\ast} \rho^2)}},
\nonumber\\
g^{R}(x, \lambda; t) &=
\frac{\Theta^{R}(\lambda/2\rho, \rho x, 2 \pi i t \rho^2)}
{\sqrt{\Theta^{D}(\lambda/2 \rho, 0, 2 \pi i t_{\ast} \rho^2)}},
\quad R=B, C, D,
\label{eqn:gR}
\end{align*}
for $t \in (0, t_{\ast})$. 
Then provided $\lambda, \lambda' \in (0, \rho), \rho>0$,
we can prove the following biorthogonality relations,
\[
\int_{\R} \overline{g^{R}(x, \lambda; t_{\ast}-t)}
g^{R}(x, \lambda'; t) dx = \delta(\lambda-\lambda'),
\quad t \in (0, t_{\ast}), 
\quad R=A, B, C, D,
\]
and the correlation kernels given in 
Lemma \ref{thm:scaling_limit_K} are written as
\begin{equation}
\cK^{R}_t(x, y; t_{\ast}, \rho)
=\int_0^{\rho} g^{R}(x, \lambda; t)
\overline{g^{R}(y, \lambda; t_{\ast}-t)} d \lambda,
\quad R=A, B, C, D.
\label{eqn:projection}
\end{equation}
Hence the four kinds of correlation kernels
obtained in the scaling limits are all
{\it reproducing kernels} with respect to the
Lebesgue measure in $\R$ for $R=A$
and in $[0, \infty)$ for $R=B, C, D$.
The formula (\ref{eqn:projection}) also implies that
$\cK^{R}_t$, $R=A, B, C, D$ can be regarded as 
{\it projection kernels}.
We can prove at least at the middle time $t=t_{\ast}/2$,
the kernels $\cK^{R}_{t_{\ast}/2}, R=A, B, C, D$ define
{\it orthogonal projections} and hence
they indeed provide
correlation kernels of determinantal point processes
\cite{Sos00,ST03,ST03b}.
More detail, see \cite{KS}.
\vskip 0.3cm

Put $t=t_{\ast}/2$ in (\ref{eqn:K_inf_A})--(\ref{eqn:K_inf_D}).
By (\ref{eqn:theta_asym}), we see that
\[
\lim_{t_{\ast} \to \infty}
\frac{\vartheta_s (\rho x+\pi i t_{\ast} \rho \lambda/2; \pi i t_{\ast} \rho^2)
\vartheta_s (\rho y - \pi i t_{\ast} \rho \lambda/2; \pi i t_{\ast} \rho^2)}
{\vartheta_2 (\pi i t_{\ast} \rho \lambda; 2  \pi i t_{\ast} \rho^2)} 
= \begin{cases}
e^{-\pi i \rho (x-y)}, \quad & \mbox{if $\lambda>0$},
\cr
e^{\pi i \rho (x-y)}, \quad & \mbox{if $\lambda<0$},
\end{cases}
\]
for $s=1,2$.
Then we obtain the following three types of limits,
\begin{align}
\cK^{A}(x, y; \rho)
&= \lim_{t_{\ast} \to \infty} \cK^{A}_{t_{\ast}/2}(x, y; t_{\ast}, \rho)
= e^{-i \pi \rho(x-y)} \int_0^{\rho} e^{2 \pi i(x-y) \lambda} d \lambda
\nonumber\\
&=\frac{\sin\{\pi \rho(x-y)\}}{\pi (x-y)},
\quad x, y \in \R, 
\label{eqn:sineA}
\\
\cK^{C}(x, y; \rho)
&= \lim_{t_{\ast} \to \infty} \cK^{R}_{t_{\ast}/2}(x, y; t_{\ast}, \rho)
\nonumber\\
&= \frac{\sin \{\pi \rho(x-y)\}}{\pi (x-y)} - \frac{\sin \{\pi \rho(x+y)\}}{\pi (x+y)},
\quad \mbox{for $R=B, C$},
\quad x, y \in [0, \infty), 
\label{eqn:sineC}
\\
\cK^{D}(x, y; \rho)
&= \lim_{t_{\ast} \to \infty} \cK^{D}_{t_{\ast}/2}(x, y; t_{\ast}, \rho)
\nonumber\\
&= \frac{\sin \{\pi \rho(x-y)\}}{\pi (x-y)} + \frac{\sin \{\pi \rho(x+y)\}}{\pi (x+y)},
\quad x, y \in [0, \infty).
\label{eqn:sineD}
\end{align}
\vskip 0.3cm
\noindent{\bf Remark 5} \,
The kernel (\ref{eqn:sineA}) is known as the sine kernel with density $\rho$,
which governs the bulk scaling limit of the determinantal point process
in GUE, 
as explained in Section \ref{sec:Introduction}.
The statistical ensemble of nonnegative square roots of
eigenvalues of $M^{\dagger} M$, 
in which $\{M\}$ are $(N+\nu) \times N$ rectangular
complex matrices and the real and imaginary parts of
their entries are independently and normally distributed,
is called the {\it chiral GUE} with parameter $\nu$.
In the scaling limit associated with $N \to \infty$
called {\it hard-edge scaling limit}, the correlation kernel of
this determinantal point process is given by
\[
K^{\rm chGUE}_{\nu}(x,y)
=\frac{2 \sqrt{xy}}{x^2-y^2}
\{ J_{\nu}(2x) y J^{\prime}_{\nu}(2y)
-J_{\nu}(2y) x J^{\prime}_{\nu}(2x) \},
\quad x, y \in [0, \infty),
\]
where $J_{\nu}(z)$ is the Bessel function and
$J^{\prime}_{\nu}(z)=dJ_{\nu}(z)/dz$ 
(see \cite{For10,KT04,KT07} and references therein).
Since $J_{1/2}(z)=\sqrt{2/(\pi z)} \sin z$
and $J_{-1/2}(z)=\sqrt{2/(\pi z)} \cos z$,
we can see that
\[
K^{\rm chGUE}_{\pm 1/2}(x,y)
= \frac{\sin\{2(x-y)\}}{\pi(x-y)}
\mp \frac{\sin\{2(x+y)\}}{\pi(x+y)}.
\]
The kernels (\ref{eqn:sineC}) and (\ref{eqn:sineD}) are
the scale changes of $K^{\rm chGUE}_{\pm1/2}$.
\vskip 0.3cm

\vskip 0.3cm
\noindent{\bf Remark 6} \,
If we take the scaling limit $N \to \infty$, $r \to \infty$
with constant density (\ref{eqn:rho1}) in the four types
of correlation kernels on the trigonometric level,
(\ref{eqn:K_f_h_A2})--(\ref{eqn:K_f_h_C}), 
the three types of sine kernels (\ref{eqn:sineA})--(\ref{eqn:sineD})
are obtained.
\vskip 0.3cm

\SSC
{Realization as Systems of Noncolliding Brownian Bridges} 
\label{sec:bridges} 
\subsection{New expressions of Macdonald denominators 
by KMLGV determinants} 
\label{sec:new_exp}

Consider the one-dimensional standard BM, 
$B(t), t \in [0, \infty)$ governed by the Wiener measure
denoted by $\rP$.
The transition probability density of BM, 
starting from $x$ at time $s$ and arriving at $y$ at time $t$,
$x, y \in \R, 0 \leq s < t < \infty$, is
denoted as $\rp^{\rm BM}(s, x; t, y)$ and defined by
\[
\rP(B(t) \in dy | B(s)=x) =
\rp^{\rm BM}(s, x; t, y) dy,
\]
with
\[
\rp^{\rm BM}(s,x; t, y) = \rp^{\rm BM}(s,y; t, x)
=\frac{1}{\sqrt{2 \pi (t-s)}} e^{-(x-y)^2/\{2(t-s)\}}.
\]
By the Markov property of BM, the Chapman-Kolmogorov equation
holds,
\begin{equation}
\int_{\R} \rp^{\rm BM}(s, x; t, y) \rp^{\rm BM}(t, y; u, z) dy
=\rp^{\rm BM}(s, x; u, z),
\quad 0 \leq s < t < u < \infty, \quad x, z \in \R.
\label{eqn:CK1}
\end{equation}

For $0 \leq s < t < \infty$, define
\begin{align}
\rp^{\rm circ}(s, x; t, y)
&= \rp^{\rm circ}(s, x; t, y; r) 
\nonumber\\
&= \begin{cases}
\displaystyle{
\sum_{w \in \Z} (-1)^w \rp^{\rm BM}(s, x; t, y+2 \pi r w)}, 
& \mbox{if $N$ is even},
\cr
\displaystyle{
\sum_{w \in \Z} \rp^{\rm BM}(s, x; t, y+2 \pi r w)}, 
& \mbox{if $N$ is odd},
\end{cases}
\nonumber\\
&= \begin{cases}
\displaystyle{
\rp^{\rm BM}(s, x; t, y) \vartheta_0 (i (x-y)r/(t-s); -1/\tau(t-s))
}, 
& \mbox{if $N$ is even},
\cr
\displaystyle{
\rp^{\rm BM}(s, x; t, y) \vartheta_3 (i (x-y)r/(t-s); -1/\tau(t-s))
}, 
& \mbox{if $N$ is odd},
\end{cases}
\nonumber\\
&= \begin{cases}
\displaystyle{
\frac{1}{2 \pi r} \vartheta_2 ( \xi(x-y); \tau(t-s))
}, 
& \mbox{if $N$ is even},
\cr
& \cr
\displaystyle{
\frac{1}{2 \pi r} \vartheta_3 ( \xi(x-y); \tau(t-s))
}, 
& \mbox{if $N$ is odd},
\end{cases}
\label{eqn:pAN1}
\end{align}
$x, y \in [0, 2\pi r)$, 
where $\xi(x)$ and $\tau(t)$ are defined by (\ref{eqn:xi_tau1})
and in the last equalities Jacobi's imaginary transformations
(\ref{eqn:Jacobi_imaginary}) were used, 
and
\begin{align}
\rp^{\rm ar}(s,x;t,y)
&= \frac{1}{\pi r} 
\sum_{n \in \Z} e^{-(n-1/2)^2 (t-s)/2r^2}
\sin \left( (n-1/2)\frac{x}{r} \right) 
\sin \left( (n-1/2)\frac{y}{r} \right) 
\nonumber\\
&= \frac{1}{2 \pi r} \Big\{
\vartheta_2(\xi(x-y); \tau(t-s)) - \vartheta_2 (\xi(x+y); \tau(t-s)) 
\Big\},
\label{eqn:pBN1}
\\
\rp^{\rm aa}(s,x;t,y)
&= \frac{1}{\pi r} 
\sum_{n \in \Z} e^{-n^2 (t-s)/2r^2}
\sin \left( \frac{nx}{r} \right) \sin \left( \frac{ny}{r} \right)
\nonumber\\
&= \sum_{k \in \Z}
\Big\{ p^{\rm BM}(s, x; t, y+2 \pi r k)
-p^{\rm BM}(s, -x; t, y+2 \pi r k) \Big\}
\nonumber\\
&= \frac{1}{2 \pi r} \Big\{
\vartheta_3(\xi(x-y); \tau(t-s)) - \vartheta_3 (\xi(x+y); \tau(t-s)) 
\Big\},
\label{eqn:pCN1}
\\
\rp^{\rm rr}(s,x;t,y)
&= \frac{1}{\pi r} 
\sum_{n \in \Z} e^{-n^2 (t-s)/2r^2}
\cos \left( \frac{nx}{r} \right) \cos \left( \frac{ny}{r} \right)
\nonumber\\
&= \sum_{k \in \Z}
\Big\{ p^{\rm BM}(s, x; t, y+2 \pi r k)
+p^{\rm BM}(s, -x; t, y+2 \pi r k) \Big\}
\nonumber\\
&= \frac{1}{2 \pi r} \Big\{
\vartheta_3(\xi(x-y); \tau(t-s)) + \vartheta_3 (\xi(x+y); \tau(t-s)) 
\Big\},
\label{eqn:pDN1}
\end{align}
$x, y \in [0, \pi r]$.
By (\ref{eqn:CK1}), we can readily confirm that,
for $0 \leq s < t < u < \infty$, 
\begin{align}
& \int_0^{2 \pi r} \rp^{\rm circ}(s, x; t, y) \rp^{\rm circ}(t, y; u, z) dy
=\rp^{\rm circ}(s, x; u, z),
\quad x, z \in [0, 2 \pi r),
\nonumber\\
& \int_0^{\pi r} \rp^{\sharp}(s, x; t, y) \rp^{\sharp}(t, y; u, z) dy
=\rp^{\sharp}(s, x; u, z),
\quad x,z \in [0, \pi r],
\quad \sharp={\rm ar, aa, rr}.
\label{eqn:CK2}
\end{align}
The functions $\rp^{\rm ar}$, $\rp^{\rm aa}$, and $\rp^{\rm rr}$ can be
interpreted as the transition probability densities
of BM in an interval $[0, \pi r]$
with absorbing boundary condition at $x=0$ and
reflecting boundary condition at $x=\pi r$, 
with absorbing boundary condition
both at $x=0$ and $x=\pi r$, 
and 
with reflecting boundary condition
both at $x=0$ and $x=\pi r$, respectively.
These facts are proved by 
the expansion of transition probability density
with eigenfunctions of Laplace equation with given boundary condition,
and by the reflection principle
of BM (see, for instance, Appendices 1.5 and 1.6 in \cite{BS02}).

We introduce time dependent $N \times N$ matrices, 
$\sfp^{\sharp}(s, \x; t, \y)$, $0 \leq s < t < \infty$, with entries
\[
(\sfp^{\sharp}(s, \x; t, \y))_{jk}
=\rp^{\sharp} (s, x_j; t, y_k),
\quad \sharp={\rm circ, ar, aa, rr},
\quad j, k =1,2, \dots, N, 
\]
for $\x=(x_1, \dots, x_N) \in \R^N$,
$\y=(y_1, \dots, y_N) \in \R^N$.
We see 
\begin{equation}
\rp^{\sharp}(0, \x; t, \y)
=\rp^{\sharp}(u-t, \y; u, \x),
\quad \sharp={\rm circ, ar, aa, rr}, 
\label{eqn:inversion}
\end{equation}
for any $0 < t < u < \infty$, $\x, \y \in \R^N$.
By the Heine identity (\ref{eqn:Heine}) in Appendix \ref{sec:appendixC}, 
the Chapman-Kolmogorov equations
(\ref{eqn:CK2}) can be extended to the following determinantal versions,
for $0 \leq s < t < u < \infty$, 
\begin{align}
\int_{\W_N^{[0, 2 \pi r)}} d \y \,
\det[\sfp^{\rm circ}(s, \x; t, \y)]
\det[\sfp^{\rm circ}(t, \y; u, \z)]
&=\det [\sfp^{\rm circ}(s, \x; u, \z)],
\quad \x, \z \in \W_N^{[0, 2 \pi r)},
\nonumber\\
\int_{\W_N^{[0, \pi r]}} d \y \,
\det[\sfp^{\sharp}(s, \x; t, \y)]
\det[\sfp^{\sharp}(t, \y; u, \z)]
&=\det [\sfp^{\sharp}(s, \x; u, \z)],
\quad \x, \z \in \W_N^{[0, \pi r]},
\nonumber\\
& \qquad 
\quad \mbox{for $\sharp=$ aa, ar, rr}.
\label{eqn:CK3}
\end{align}
The determinant $\det[\sfp^{\rm circ}(s, \x; t, \y)]$
with $\x, \y \in \W_N^{[0, 2 \pi r)}$,
$0 \leq s < t $, 
is the KMLGV determinant giving the total probability mass
of $N$-tuple of noncolliding Brownian paths
on a circle with radius $r>0$, 
starting from the unlabeled configuration $\x$ at time $s$ and 
arriving at the unlabeled configuration $\y$
at time $t>s$
\cite{For90,Ful04,LW16}.
The determinants $\det[\sfp^{\rm ar}(s, \x; t, \y)]$, 
$\det[\sfp^{\rm aa}(s, \x; t, \y)]$, and
$\det[\sfp^{\rm rr}(s, \x; t, \y)]$
can be regarded as the KMLGV determinants
for the noncolliding BMs in the interval
$[0, \pi r]$
with absorbing boundary condition at $x=0$ and
reflecting boundary condition at $x=\pi r$, 
with absorbing boundary condition
both at $x=0$ and $x=\pi r$, 
and 
with reflecting boundary condition
both at $x=0$ and $x=\pi r$, respectively.
See \cite{TW07,LW17} for the noncolliding Brownian bridges
starting from and returning to the origin
with reflecting or absorbing walls. 

We consider the following seven types of configurations of $N$ points,
$\v^{\RN}=\v^{\RN}(r)$ with the elements,
\begin{align}
v^{\AN}_j &= v^{\AN}_j(r) = 
\frac{2 \pi r}{N}(j-1), 
\nonumber\\
v^{\RN}_j &=
v^{\RN}_j(r) = \frac{2 \pi r}{\cN^{\RN}}(j-1/2),
\quad \mbox{for $\RN=\BN, \BNv$},
\nonumber\\
v^{\RN}_j &= 
v^{\RN}_j(r) = \frac{2 \pi r}{\cN^{\RN}} j,
\quad \mbox{for $\RN=\CN, \CNv, \BCN$},
\nonumber\\
v^{\DN}_j &= 
v^{\DN}_j(r) = \frac{\pi r}{N-1} (j-1),
\qquad j=1,2, \dots, N.
\label{eqn:v}
\end{align}
The configurations $\v^{\RN}$ make equidistant series of points
in $[0, 2 \pi r)$ for $\RN=\AN$
and in $[0, \pi r]$ for others.
We also consider $N \times N$ matrices whose entries
are given by the biorthogonal theta functions
studied in Section \ref{sec:orthogonality},
\[
\sfM^{\RN}(\x,t)= \Big( M^{\RN}_j(x_k, t) \Big)_{1 \leq j, k \leq N}.
\]
Then the following relations hold between matrices.
\begin{lem}
\label{thm:KM_matrix}
Consider the $N \times N$ matrices
$\sfr^{\RN}(t)$ with the following entries;
for $j=1,\dots, N$, 
\begin{align*}
(\sfr^{\AN}(t))_{jk}
&= \frac{2 \pi r}{\cN^{\AN}}
e^{-\pi i(J^{\AN}(j))^2 \tau(t) - i (\cN^{\AN}-2J^{\AN}(j))v^{\AN}_k/2r }, 
\quad k=1, \dots, N, 
\nonumber\\
(\sfr^{\BN}(t))_{jk}
&=\begin{cases}
\displaystyle{
\frac{4 \pi r}{\cN^{\BN}}
e^{-\pi i (J^{\BN}(j))^2 \tau(t)}
\sin \left[ 
(\cN^{\BN}-2J^{\BN}(j)) \frac{v^{\BN}_k}{2r}
\right]
},
& k=1, \dots, N-1,
\cr
\displaystyle{
\frac{2 \pi r}{\cN^{\BN}}
e^{-\pi i (J^{\BN}(j))^2 \tau(t)}
\sin[(\cN^{\BN}-2 J^{\BN}(j))\pi/2]
},
& k=N.
\end{cases}
\nonumber\\
\nonumber\\
(\sfr^{\BNv}(t))_{jk}
&= \frac{4 \pi r}{\cN^{\BNv}}
e^{-\pi i (J^{\BNv}(j))^2 \tau(t)}
\sin \left[ 
(\cN^{\BNv}-2J^{\BNv}(j)) \frac{v^{\BNv}_k}{2r}
\right],
\quad k=1,\dots, N, 
\nonumber\\
(\sfr^{\RN}(t))_{jk}
&= \frac{4 \pi r}{i \cN^{\RN}}
e^{-\pi i (J^{\RN}(j))^2 \tau(t)}
\sin \left[ 
(\cN^{\RN}-2J^{\RN}(j)) \frac{v^{\RN}_k}{2r}
\right],
\quad k=1,\dots, N, 
\nonumber\\
& \hskip 7cm
\mbox{for $\RN=\CN, \BCN$},
\end{align*}
\begin{align}
(\sfr^{\CNv}(t))_{jk}
&=\begin{cases}
\displaystyle{
\frac{4 \pi r}{i \cN^{\CNv}}
e^{-\pi i (J^{\CNv}(j))^2 \tau(t)}
\sin \left[ 
(\cN^{\CNv}-2J^{\CNv}(j)) \frac{v^{\CNv}_k}{2r}
\right]
},
& k=1, \dots, N-1,
\cr
\displaystyle{
\frac{2 \pi r}{i \cN^{\CNv}}
e^{-\pi i (J^{\CNv}(j))^2 \tau(t)}
\sin[(\cN^{\CNv}-2J^{\CNv}(j)) \pi/2)]
},
& k=N.
\end{cases}
\nonumber\\
(\sfr^{\DN}(t))_{jk}
&=\begin{cases}
\displaystyle{ \frac{2 \pi r}{\cN^{\DN}} e^{-\pi i (J^{\DN}(j))^2 \tau(t)}},
& k=1,
\cr
\displaystyle{\frac{4 \pi r}{\cN^{\DN}} e^{- \pi i (J^{\DN}(j))^2 \tau(t)} 
\cos \left[(\cN^{\DN}-2J^{\DN}(j)) \frac{v^{\DN}_k}{2r}  \right]},
& k=2, \dots, N-1,
\cr
\displaystyle{
\frac{2 \pi r}{\cN^{\DN}} e^{-\pi i (J^{\DN}(j))^2 \tau(t)} 
\cos[(\cN^{\DN}-2 J^{\DN}(j))\pi/2] },
& k=N.
\end{cases}
\label{eqn:r}
\end{align}
Then for $t \in [0, \infty)$, 
\begin{align}
\sfr^{\AN}(t) \sfp^{\rm circ}(0, \v^{\AN}; t, \x) &= \sfM^{\AN}(\x, t),
\quad \x \in \W_N^{[0, 2 \pi r)},
\label{eqn:matrix_relation_A}
\\
\sfr^{\RN}(t) \sfp^{\rm ar}(0, \v^{\RN}; t, \x) &= \sfM^{\RN}(\x, t),
\quad \x \in \W_N^{[0, \pi r]},
\quad \RN=\BN, \CNv, \BCN, 
\label{eqn:matrix_relation_absref}
\\
\sfr^{\RN}(t) \sfp^{\rm aa}(0, \v^{\RN}; t, \x) &= \sfM^{\RN}(\x, t),
\quad \x \in \W_N^{[0, \pi r]},
\quad \RN=\BNv, \CN, 
\label{eqn:matrix_relation_abs}
\\
\sfr^{\DN}(t) \sfp^{\rm rr}(0, \v^{\DN}; t, \x) &= \sfM^{\DN}(\x, t),
\quad \x \in \W_N^{[0, \pi r]}.
\label{eqn:matrix_relation_ref}
\end{align}
\end{lem}
\vskip 0.3cm
\vskip 0.3cm
\noindent{\bf Remark 7} \,
Before proving the equalities 
(\ref{eqn:matrix_relation_A})--(\ref{eqn:matrix_relation_ref}), 
we explain how one of
the boundary conditions, `ar', `aa', `rr', is chosen
for each $\RN=\BN$, $\BNv$, $\CN$, $\CNv$, $\BCN$, $\DN$.
As given by the Macdonald denominator formulas 
(\ref{eqn:Macdonald3}), 
the determinants of the matrices in RHS of 
(\ref{eqn:matrix_relation_A})--(\ref{eqn:matrix_relation_ref})
are proportional to the Macdonald denominators 
(\ref{eqn:Macdonald_denominators}) with (\ref{eqn:xi_tau1}). 
If $\RN=\BNv$, $\CN$, then due to the factors,
$\prod_{\ell=1}^{N} \vartheta_1(2 \xi_{\ell}(\x); \cdot)=
\prod_{\ell=1}^{N} \vartheta_1(x_{\ell}/\pi r; \cdot)$,
$W^{\RN}(\xi(\x); \tau)=0$ when $x_{\ell}=0$ or
$x_{\ell}=\pi r$ for any $\ell \in \{1,2, \dots, N\}$.
While, if $\RN=\BN, \CNv, \BCN$, 
then due to the factors,
$\prod_{\ell=1}^{N} \vartheta_1(\xi_{\ell}(\x); \cdot)=
\prod_{\ell=1}^{N} \vartheta_1(x_{\ell}/2\pi r; \cdot)$,
$W^{\RN}(\xi(\x); \tau)=0$ when $x_{\ell}=0$ 
for any $\ell \in \{1,2, \dots, N\}$, but 
$W^{\RN}(\xi(\x); \tau) \not=0$ when $x_{\ell} = \pi r,
\ell \in \{1,2, \dots, N\}$. 
There is no such factor in $W^{\DN}(\xi(\x); \tau)$.
Therefore, the noncolliding BMs in an interval $[0, \pi r]$,
whose KMLGV determinants are
proportional to $\det[\sfM^{\RN}(\x, t)]$
(as (\ref{eqn:det_relation}) given below), 
could be considered under the boundary condition
`aa' for $\RN=\BNv, \CN$,
`ar' for $\RN=\BN, \CNv, \BCN$, and
`rr' for $\RN=\DN$, respectively.
Moreover, we see in the initial configurations
$\v^{\RN}$ at time $t=0$ given by (\ref{eqn:v}) that
$v^{\RN}_1 > 0$ for $\RN=\BN$, $\BNv$, $\CN$, $\CNv$, $\BCN$,
while $v^{\DN}_1=0$,
and that
$v^{\RN}_N < \pi r$ for $\RN=\BNv, \CN, \BCN$, while
$v^{\RN}_N=\pi r$ for $\RN=\BN, \CNv, \DN$.
They are consistent with the choice of boundary conditions.
\vskip 0.3cm
\noindent{\it Proof of Lemma \ref{thm:KM_matrix}} \
First we prove (\ref{eqn:matrix_relation_A}) when $N$ is even.
By the definitions of 
$\sfr^{\AN}(t)$ and $\sfp^{\rm circ}(0, \v^{\AN}; t, \x)$
with (\ref{eqn:v}),
the $(j,k)$-entry of 
LHS of (\ref{eqn:matrix_relation_A}) is given by
\begin{align*}
L^{\AN}_{jk} &= 
\frac{2 \pi r}{N} e^{-\pi i (j-1/2)^2 \tau(t)} 
\nonumber\\
& \quad
\times \sum_{\ell=1}^N e^{2 \pi i(j-1/2-N/2)(\ell-1)/N}
\frac{1}{2 \pi r} 
\vartheta_2 \left( \frac{\ell-1}{N}-\frac{x_k}{2 \pi r}; \tau(t) \right)
\nonumber\\
&=
\frac{1}{N} e^{-\pi i (j-1/2)^2 \tau(t)}
\sum_{n \in \Z} 
e^{ (n-1/2)^2 \tau(t) \pi i - (2n-1) \{x_k/(2 \pi r)\} \pi i}
\nonumber\\
& \quad
\times \sum_{\ell=1}^N e^{2 \pi i (\ell-1) \{(j-1/2-N/2)+(n-1/2) \}/N},
\quad 1 \leq j, k \leq N, 
\end{align*}
where we have used the definition of $\vartheta_2$
given by (\ref{eqn:theta}).
We note $(j-1/2-N/2)+(n-1/2) \in \Z$
for $N$ even and use the equality
\begin{equation}
\sum_{\ell=1}^N e^{2 \pi i (\ell-1) \theta/N}
= N \sum_{m \in \Z} {\bf 1}(\theta+Nm=0),
\quad \theta \in \Z, \quad N \in \N.
\label{eqn:sum_equality_2}
\end{equation}
Then we obtain
\[
L^{\AN}_{jk} =
e^{- \pi i (j-1/2)^2 \tau(t)}
\sum_{m \in \Z}
e^{
\{ m N + (j-1/2-N/2) \}^2 \tau(t) \pi i 
+2 \{ mN + (j-1/2-N/2) \}
\{x_{k}/(2 \pi r) \} \pi i}.
\]
It is easy to confirm the equality
\begin{align*}
& -\pi i (j-1/2)^2 \tau(t) 
+\{mN+(j-1/2-N/2) \}^2 \tau(t) \pi i
+2 \{mN+(j-1/2-N/2)\} \frac{x_{k}}{2 \pi r} \pi i
\nonumber\\
&= 
(m-1/2)^2 N^2 \tau(t) \pi i
+ (2m-1) N
\left\{ (j-1/2) \tau(t)+ \frac{x_k}{2 \pi r} \right\} \pi i
+ 2 \pi i (j-1/2) \frac{x_k}{2 \pi r}.
\end{align*}
Hence we have the equality
$L_{jk}=M^{\AN}_j(x_k, t)$.
We can similarly prove (\ref{eqn:matrix_relation_A}) for odd $N$. \\
Next we explain how to prove (\ref{eqn:matrix_relation_abs})
for $\RN=\CN$. 
The $(j,k)$-entry of LHS is
\begin{align*}
L^{\CN}_{jk} &=
\frac{i e^{-\pi i j^2 \tau(t)} }{N+1} \left[
\sum_{\ell=1}^N \sin \left[ \frac{\pi \{j-(N+1) \} \ell}{N+1} \right]
\vartheta_3 \left( \frac{\ell}{2 (N+1)} - \frac{x_k}{2 \pi r} ; \tau(t) \right)
\right.
\nonumber\\
& \qquad \qquad \qquad
\left.
-\sum_{\ell=1}^N 
\sin \left[ \frac{\pi \{j-(N+1) \} \ell}{N+1} \right]
\vartheta_3 \left( \frac{\ell}{2 (N+1)} + \frac{x_k}{2 \pi r} ; \tau(t) \right)
\right], \quad 1 \leq j, k \leq N.
\end{align*}
By the fact that $\sin[\pi\{j-(N+1)\} \ell/(N+1)]=0$ when
$\ell=0$ and $\ell=N+1$, 
and by the parity $\sin(-\pi v)=-\sin(\pi v)$,
$\vartheta_3(-v; \tau)=\vartheta_3(v; \tau)$, 
this entry is equal to
\begin{align*}
L^{\CN}_{jk} &=
\frac{i e^{-\pi i j^2 \tau(t)} }{N+1} \left[
\sum_{\ell=0}^{N+1} \sin \left[ \frac{\pi \{j-(N+1) \} \ell}{N+1} \right]
\vartheta_3 \left( \frac{\ell}{2 (N+1)} - \frac{x_k}{2 \pi r} ; \tau(t) \right)
\right.
\nonumber\\
& \qquad \qquad \qquad
\left.
-\sum_{\ell=1}^N 
\sin \left[ \frac{\pi \{j-(N+1) \} \ell}{N+1} \right]
\vartheta_3 \left( \frac{\ell}{2 (N+1)} + \frac{x_k}{2 \pi r} ; \tau(t) \right)
\right]
\nonumber\\
&=
\frac{i e^{-\pi i j^2 \tau(t)} }{N+1} 
\sum_{\ell=-N}^{N+1} \sin \left[ \frac{\pi (j -N-1) \ell}{N+1} \right]
\vartheta_3 \left( \frac{\ell}{2 (N+1)} - \frac{x_k}{2 \pi r} ; \tau(t) \right),
\quad 1 \leq j,k \leq N.
\end{align*}
Note that, in the last expression, the number of terms of
summation is equal to $2N+2=\cN^{\CN}$. 
Then we use the definition (\ref{eqn:theta}) of the Jacobi theta function
$\vartheta_3$, and rewrite the above as
\begin{align*}
L^{\CN}_{jk} &= \frac{e^{-\pi i j^2 \tau(t)} }{2 (N+1)}
\sum_{n \in \Z} e^{\{n^2 \tau(t) -2 n x_{k}/(2 \pi r)\} \pi i}
\nonumber\\
& \quad \times
\sum_{\ell=-N}^{N+1}
\Big\{
e^{2 \pi i \ell \{ j-(N+1)+n\}/\{2(N+1)\} }
- e^{- 2 \pi i \ell \{j-(N+1)-n\}/\{2(N+1)\}} \Big\},
\quad 1 \leq j, k \leq N.
\end{align*}
By the equality (\ref{eqn:sum_equality_2}) with the replacement 
$N \to \cN^{\CN}=2(N+1)$, 
we can verify 
$L_{jk}=M^{\CN}_j(x_k, t)$.
For other types of $\RN$ than $\AN$ and $\CN$, 
we can show that 
the $(j,k)$-entries of LHS of 
(\ref{eqn:matrix_relation_absref}), (\ref{eqn:matrix_relation_abs}),
and (\ref{eqn:matrix_relation_ref}) are written as
\begin{align*}
L^{\BN}_{jk} &=
-\frac{e^{-\pi i (j-1)^2 \tau(t)}}{2N-1}
\sum_{\ell=-N+2}^N
\sin\left[ \frac{2\pi(j-N-1/2)(\ell-1/2)}{2N-1} \right]
\nonumber\\
& \qquad \qquad \times
\left\{ \vartheta_2 \left(
\frac{\ell-1/2}{2N-1}-\frac{x_k}{2 \pi r}; \tau(t) \right)
- \vartheta_2 \left(
\frac{\ell-1/2}{2N-1}+\frac{x_k}{2 \pi r}; \tau(t) \right) \right\},
\nonumber\\
L^{\BNv}_{jk} &=
-\frac{e^{-\pi i (j-1)^2 \tau(t)}}{N}
\sum_{\ell=-N+1}^N
\sin\left[ \frac{\pi(j-N-1)(\ell-1/2)}{N} \right]
\vartheta_3 \left(
\frac{\ell-1/2}{2N}-\frac{x_k}{2 \pi r}; \tau(t) \right),
\nonumber\\
L^{\CNv}_{jk} &=
\frac{i e^{-\pi i (j-1/2)^2 \tau(t)}}{2N}
\sum_{\ell=-N+1}^N
\sin\left[ \frac{\pi(j-N-1/2)\ell}{N} \right]
\nonumber\\
& \qquad \qquad \times
\left\{ \vartheta_2 \left(
\frac{\ell}{2N}-\frac{x_k}{2 \pi r}; \tau(t) \right)
- \vartheta_2 \left(
\frac{\ell}{2N}+\frac{x_k}{2 \pi r}; \tau(t) \right) \right\},
\nonumber\\
L^{\BCN}_{jk} &=
\frac{2 i e^{-\pi i j^2 \tau(t)}}{2N+1}
\sum_{\ell=-N}^N
\sin\left[ \frac{2 \pi(j-N-1/2) \ell}{2N+1} \right]
\vartheta_2 \left(
\frac{\ell}{2N+1}-\frac{x_k}{2 \pi r}; \tau(t) \right),
\nonumber\\
L^{\DN}_{jk} &=
\frac{e^{-\pi i (j-1)^2 \tau(t)}}{2(N-1)}
\sum_{\ell=-N+2}^{N-1}
\cos\left[ \frac{\pi(j-N)\ell}{N-1} \right]
\nonumber\\
& \qquad \qquad \times
\left\{ \vartheta_3 \left(
\frac{\ell}{2(N-1)}-\frac{x_k}{2 \pi r}; \tau(t) \right)
+ \vartheta_3 \left(
\frac{\ell}{2(N-1)}+\frac{x_k}{2 \pi r}; \tau(t) \right) \right\}.
\end{align*}
For each type of $\RN$, the number of terms in the summation
is equal to $\cN^{\RN}$.
Hence we can apply the equality (\ref{eqn:sum_equality_2}) 
with the replacement 
$N \to \cN^{\RN}$. 
In this way we can  prove
(\ref{eqn:matrix_relation_absref})
for $\RN=\BN, \CNv, \BCN$, 
(\ref{eqn:matrix_relation_abs})
also for $\RN=\BNv$,
and (\ref{eqn:matrix_relation_ref}). 
The proof is complete.
\qed
\vskip 0.3cm

If we take the determinants of both sides
of the equalities 
(\ref{eqn:matrix_relation_A})--(\ref{eqn:matrix_relation_ref}), 
we obtain the equalities
\begin{equation}
\det[\sfr^{\RN}(t)]
\det[\sfp^{\sharp}(0, \v^{\RN}; t, \x)]
=\det[\sfM^{\RN}(\x, t)],
\label{eqn:det_relation}
\end{equation}
where $\sharp=$ circ for $\RN=\AN$,
$\sharp=$ ar for $\RN=\BN, \CNv, \BCN$,
$\sharp=$ aa for $\RN=\BNv, \CN$,
and $\sharp=$ rr for $\RN=\DN$.
Combine them with the Macdonald denominator
formulas of Rosengren and Schlosser \cite{RS06},
which are written as (\ref{eqn:Macdonald3})
in the present paper, 
we obtain new determinantal expressions
for the Macdonald denominators.

\begin{prop}
\label{thm:KM_det0}
For the irreducible reduced affine root systems,
$\RN=\AN$, $\BN$, $\BNv$, $\CN$, $\CNv$, $\BCN$, $\DN$, 
the Macdonald denominators
$W^{\RN}$ defined by
{\rm (\ref{eqn:Macdonald_denominators})}
are proportional to the KMLGV determinants 
for noncolliding Brownian paths
starting from the configurations
$\v^{\RN}$ given by {\rm (\ref{eqn:v})} as follows.
Let $s(N)=0$ if $N$ is even, and 
$s(N)=3$ if $N$ is odd, then
\begin{align}
&\vartheta_{s(N)} \left( \sum_{j=1}^N \xi(x_j); N \tau(t) \right)
W^{\AN}(\xi(\x); N \tau(t))
= b^{\AN}(t) \det[ \sfp^{\rm circ}(0, \v^{\AN}; t, \x)],
\label{eqn:new_Mac_A}
\\
&W^{\RN}(\xi(\x); \cN^{\RN} \tau(t))
=\begin{cases}
b^{\RN}(t) \det[ \sfp^{\rm ar}(0, \v^{\RN}; t, \x)], \quad
\mbox{for $\RN=\BN, \CNv, \BCN$},
\cr
b^{\RN}(t) \det[ \sfp^{\rm aa}(0, \v^{\RN}; t, \x)], \quad
\mbox{for $\RN=\BNv, \CN$},
\cr
b^{\DN}(t) \det[ \sfp^{\rm rr}(0, \v^{\DN}; t, \x)], \quad
\mbox{for $\RN=\DN$},
&
\end{cases}
\label{eqn:new_Mac_RN}
\end{align}
with the coefficients
\[
b^{\RN}(t)= \begin{cases}
\displaystyle{
 i^{N(N+1)/2} 
\frac{\det[\sfr^{\AN}(t)]}{a^{\AN}(t)} },
& \mbox{for $\RN=\AN$, if $N$ is even},
\cr
\displaystyle{
 i^{(N-1)(N-2)/2} 
\frac{\det[\sfr^{\AN}(t)]}{a^{\AN}(t)} },
& \mbox{for $\RN=\AN$, if $N$ is odd},
\cr
\displaystyle{
\frac{\det[\sfr^{\RN}(t)]}{a^{\AN}(t)}},
& \mbox{for $\RN=\BN, \BNv, \DN$},
\cr
\displaystyle{
i^N
\frac{\det[\sfr^{\RN}(t)]}{a^{\AN}(t)}},
& \mbox{for $\RN=\CN,\CNv,\BCN$},
\end{cases}
\]
where the factors $a^{\RN}(t)$ and 
the entries of the matrices $\sfr^{\RN}(t)$
are given by (\ref{eqn:a1}) and (\ref{eqn:r}),
respectively.
\end{prop}
\vskip 0.3cm
\noindent{\bf Remark 8} \,
Forrester proved the equality (\ref{eqn:new_Mac_A})
for the type $\AN$ independently of the
Macdonald denominator formulas given by
Rosengren and Schlosser \cite{RS06}.
For $N$ even (resp. $N$ odd),
(\ref{eqn:new_Mac_A}) is a special case with
$\alpha=1/2+1/N$
(resp. $\alpha=1/N$) of Eq.(5.111)
(resp. Eq.(5.110))
in Proposition 5.6.3 in \cite{For10}.
The matrix relation (\ref{eqn:matrix_relation_A})
was also used to prove (\ref{eqn:det_relation})
for $\RN=\AN$ in pages 216-217 in \cite{For10}.
Moreover, explicit evaluation of 
$\det[\sfr^{\AN}(t)]$ was found there.
If we use this result, we obtain
\begin{align*}
b^{\AN}(t)
&= \frac{(2 \pi r)^N}{N^{N/2}}
q(N \tau(t))^{(N-1)(N-2)/24}
q_0(N \tau(t))^{(N-1)(N-2)/2}
\nonumber\\
&= \frac{(2 \pi r)^N}{N^{N/2}}
\eta(N \tau(t))^{(N-1)(N-2)/2},
\end{align*}
where $\eta(\tau)$ is the Dedekind modular function
(see, for instance, Sec.23.15 in \cite{NIST10}),
\[
\eta(\tau)=
q(\tau)^{1/12} q_0(\tau)
=e^{\tau \pi i/12} \prod_{n=1}^{\infty} (1-e^{2n \tau \pi i}).
\]
See also \cite{For06} and references therein. 
Lemma \ref{thm:KM_matrix}
and Proposition \ref{thm:KM_det0}
are extensions of Forrester's results
to other six types of matrices and their determinants.
Here we identify LHS of the equations
(\ref{eqn:new_Mac_A}) and (\ref{eqn:new_Mac_RN})
as the Macdonald denominators
and the determinants in RHS of them
as the KMLGV determinants of
noncolliding Brownian paths.
\vskip 0.3cm

\subsection{Noncolliding Brownian bridges} 
\label{sec:nB_bridges}

The following is derived
by Lemma \ref{thm:KM_matrix}.

\begin{prop}
\label{thm:KM_det}
The probability densities given by {\rm (\ref{eqn:PtR1})}
for the determinantal point processes,
$(\Xi^{\RN}, \bP^{\RN}_t, t \in (0, t_{\ast}))$,
have the following expressions,
\begin{align}
\bp^{\AN}_t(\x)
& =\frac{
\det[\sfp^{\rm circ}(0, \v^{\AN}; t, \x)]
\det[\sfp^{\rm circ}(t, \x; t_{\ast}, \v^{\AN})]}
{\det[\sfp^{\rm circ}(0, \v^{\AN}; t_{\ast}, \v^{\AN})]},
\quad \x \in \W_N^{[0, 2 \pi r)},
\nonumber\\
\bp^{\RN}_t(\x)
&=\frac{
\det[\sfp^{\rm ar}(0, \v^{\RN}; t, \x)]
\det[\sfp^{\rm ar}(t, \x; t_{\ast}, \v^{\RN})]}
{\det[\sfp^{\rm ar}(0, \v^{\RN}; t_{\ast}, \v^{\RN})]},
\quad \x \in \W_N^{[0, \pi r]},
\quad \mbox{for $\RN=\BN,\CNv,\BCN$},
\nonumber\\
\bp^{\RN}_t(\x)
&=\frac{
\det[\sfp^{\rm aa}(0, \v^{\RN}; t, \x)]
\det[\sfp^{\rm aa}(t, \x; t_{\ast}, \v^{\RN})]}
{\det[\sfp^{\rm aa}(0, \v^{\RN}; t_{\ast}, \v^{\RN})]},
\quad \x \in \W_N^{[0, \pi r]},
\quad \mbox{for $\RN=\BNv,\CN$},
\nonumber\\
\bp^{\DN}_t(\x)
&=\frac{
\det[\sfp^{\rm rr}(0, \v^{\DN}; t, \x)]
\det[\sfp^{\rm rr}(t, \x; t_{\ast}, \v^{\DN})]}
{\det[\sfp^{\rm rr}(0, \v^{\DN}; t_{\ast}, \v^{\DN})]},
\quad \x \in \W_N^{[0, \pi r]}.
\label{eqn:PtR2}
\end{align}
\end{prop}
\noindent{\it Proof} \,
In the equalities (\ref{eqn:det_relation}), 
if we replace $t$ by $t_{\ast}-t$ 
and consider the
complex conjugate of the obtained equalities,
then by (\ref{eqn:inversion}), we have
\[
\det[\overline{\sfr^{\RN}(t_{\ast}-t)}]
\det[\sfp^{\sharp}(t, \x; t_{\ast}, \v^{\RN})]
=\det[\overline{\sfM^{\RN}(\x, t_{\ast}-t)}].
\]
Hence (\ref{eqn:PtR1}) with (\ref{eqn:q1}) gives
\[
\bp^{\RN}_t(\x)
= c^{\RN} \det[ \sfp^{\sharp}(0, \v^{\RN}; t, \x)]
\det[ \sfp^{\sharp}(t, \x; t_{\ast}, \v^{\RN})]
\]
with constants $c^{\RN}$ which do not depend on $\x$.
Since $\bp^{\RN}_t(x)$ is normalized 
as (\ref{eqn:normalization_p}), 
the Chapman-Kolmogorov equations
(\ref{eqn:CK3}) determine the constants as
$c^{\RN}=1/\det[\sfp^{\sharp}(0, \v^{\RN}; t_{\ast}, \v^{\RN})]$.
The proof is complete.
\qed
\vskip 0.3cm

From the expressions (\ref{eqn:PtR2}) in 
Proposition \ref{thm:KM_det}, 
we can conclude the following.
(See, for instance, Part I, IV.4.22 of \cite{BS02} and
Section V.C of \cite{KT04} 
for Brownian bridges.) 

\begin{thm}
\label{thm:B_bridges} 
{\rm (i)} 
The one-parameter family of determinantal point process, 
$(\Xi^{\AN} \bP^{\AN}_t, t \in (0, t_{\ast}))$, 
is realized as the particle configuration at 
time $t \in (0, t_{\ast})$ of the noncolliding Brownian 
bridges on a circle with radius $r$,
starting from and returning to
the configuration $\v^{\AN}=(2 \pi r (j-1)/N)_{j=1}^N$. 
\\
{\rm (ii)} 
For $\RN=\BN, \CNv, \BCN$, 
each one-parameter family of determinantal point process, 
$(\Xi^{\RN}, \bP^{\RN}_t, t \in (0, t_{\ast}))$, 
is realized as the particle configuration at 
time $t \in (0, t_{\ast})$ of the noncolliding Brownian 
bridges starting from and returning to
the configuration $\v^{\RN}$
given by {\rm (\ref{eqn:v})} in the interval
$[0, \pi r]$ with absorbing boundary condition
at $x=0$ and reflecting boundary condition at $x=\pi r$.
\\
{\rm (iii)} 
For $\RN=\BNv, \CN$, 
each one-parameter family of determinantal point process, 
$(\Xi^{\RN}, \bP^{\RN}_t, t \in (0, t_{\ast}))$, 
is realized as the particle configuration at 
time $t \in (0, t_{\ast})$ of the noncolliding Brownian 
bridges starting from and returning to
the configuration $\v^{\RN}$
given by {\rm (\ref{eqn:v})} in the interval
$[0, \pi r]$ with absorbing boundary condition
at both edges. \\
{\rm (iv)} The one-parameter family of determinantal point process, 
$(\Xi^{\DN}, \bP^{\DN}_t, t \in (0, t_{\ast}))$, 
is realized as the particle configuration at 
time $t \in (0, t_{\ast})$ of the noncolliding 
bridges starting from and returning to
the configuration $\v^{\DN}=(\pi r(j-1)/(N-1))_{j=1}^N$
in the interval
$[0, \pi r]$ with the reflecting boundary conditions
at both edges. 
\end{thm}

\SSC
{Concluding Remarks} 
\label{sec:remarks} 

In the present paper we have constructed
seven types of one-parameter 
families of determinantal point processes,
$(\Xi^{\RN}, \bP^{\RN}_t, t \in (0, t_{\ast}))$,
$\RN=\AN$, $\BN$, $\BNv$, $\CN$, $\CNv$, $\BCN$, $\DN$.
These point processes can be interpreted as configurations
at time $t \in (0, t_{\ast})$ of the noncolliding Brownian
bridges starting from and returning to the equidistant configurations
$\v^{\RN}$ given by (\ref{eqn:v}).
In this picture, the variety of elliptic determinantal processes
is due to various choices of configurations pinned
at the initial time $t=0$ and at the final time $t=t_{\ast}$.
If we regard these Brownian bridges
on a circle with radius $r$, $\rP^1(r)$, or in an interval $[0, \pi r]$
with time duration $t_{\ast}$ as the statistical ensembles
of noncolliding paths on the spatio-temporal 
cylinder $\rP^1(r) \times (0, t_{\ast})$ 
or on the spatio-temporal plane
$[0, \pi r] \times (0, t_{\ast})$,
$\v^{\RN}$ gives a boundary condition to the paths.
The degeneracy of types  in the scaling limit
$N \to \infty$, $r \to \infty$ with constant 
density $\rho$ of paths shown by Theorem \ref{thm:infinite_systems}
is caused by vanishing of the boundary effect
in this bulk limit. 

In previous papers \cite{Kat15,Kat16,Kat17},
the processes associated with the affine root systems
$\AN$, $\BN$, $\CN$, and $\DN$ were characterized as
solutions of some systems of stochastic differential
equations (SDEs).
Characterization of the present determinantal
point processes $(\Xi^{\RN}, \bP^{\RN}_t, t \in (0, t_{\ast}))$
in terms of SDEs should be further studied.
The noncolliding Brownian bridges discussed in
Section \ref{sec:bridges} are
determinantal \cite{BR05,KT07,Kat15_Springer}, 
and thus the spatio-temporal correlation kernels
should be determined. 

As mentioned in Section \ref{sec:Introduction} 
and in Remark 5 in Section \ref{sec:infinite},
the present determinantal point processes
are elliptic extensions of the eigenvalue ensembles
of Hermitian random matrices in GUE and chiral GUE.
The trigonometric reductions discussed 
in Section \ref{sec:temp_homo_finite} are
related with the eigenvalue distributions of
random matrices in ${\rm U}(N)$ (CUE)
\cite{Meh04,For10}
or in other orthogonal and symplectic matrices
(see Remark 3 in Section \ref{sec:temp_homo_finite} 
and Section 2.3 (c) in \cite{Sos00}). 
It is an interesting future problem
to find the statistical ensembles of random matrices
on the elliptic level whose eigenvalues 
realize the present seven types of elliptic 
determinantal point processes.

In \cite{For06} Forrester studied the quantum $N$-particle
systems in two dimensions with doubly periodic boundary conditions,
in which the $N$-body potentials and wave functions
are described using the Jacobi theta functions.
He constructed the doubly periodic probability measures
on a complex plane and discussed 
solvability and universality of the obtained
two-dimensional systems.
From the view point of the present study,
his systems are of type $\AN$ and
they are truly elliptic.
Generalization of his study to the two-dimensional systems
associated with the other six types
of irreducible reduced affine root systems
is reported in \cite{Kat18b}. 

\vskip 0.5cm
\noindent{\bf Acknowledgements} \,
On sabbatical leave from Chuo University, 
this study was done in 
Fakult\"{a}t f\"{u}r Mathematik, Universit\"{a}t Wien,
in which the present author thanks Christian Krattenthaler very much
for his hospitality. 
The author also thanks Michael Schlosser and
Peter John Forrester for useful comments on the
manuscript.
He expresses his gratitude to Piotr Graczyk and Jacek Ma{\l}ecki 
for valuable discussion on noncolliding Brownian bridges,
which enabled him to give correct statements 
in Section \ref{sec:bridges}.
This work was supported by
the Grant-in-Aid for Scientific Research (C) (No.26400405),
(B) (No.18H01124), and 
(S) (No.16H06338) 
of Japan Society for the Promotion of Science.

\appendix
\SSC{The Jacobi Theta Functions}
\label{sec:appendixA}
Let
\[
z=e^{v \pi i}, \quad q=e^{\tau \pi i},
\]
where $v, \tau \in \C$ and $\Im \tau > 0$. 
The Jacobi theta functions are defined as follows \cite{WW27,NIST10}, 
\begin{align}
\vartheta_0(v; \tau) &= 
-i e^{\pi i (v+\tau/4)} \vartheta_1 \left( v + \frac{\tau}{2}; \tau \right)=
\sum_{n \in \Z} (-1)^n q^{n^2} z^{2n} 
=1+ 2 \sum_{n=1}^{\infty}(-1)^n e^{\tau \pi i n^2} \cos(2 n \pi v),
\nonumber\\
\vartheta_1(v; \tau) &= i \sum_{n \in \Z} (-1)^n q^{(n-1/2)^2} z^{2n-1}
=2 \sum_{n=1}^{\infty} (-1)^{n-1} e^{\tau \pi i (n-1/2)^2} \sin\{(2n-1) \pi v\},
\nonumber\\
\vartheta_2(v; \tau) 
&= \vartheta_1 \left( v+ \frac{1}{2}; \tau \right)=
\sum_{n \in \Z} q^{(n-1/2)^2} z^{2n-1}
=2 \sum_{n=1}^{\infty} e^{\tau \pi i (n-1/2)^2} \cos \{(2n-1) \pi v\},
\nonumber\\
\vartheta_3(v; \tau) 
&= 
e^{\pi i (v+\tau/4)} \vartheta_1 \left( v+\frac{1+\tau}{2}; \tau \right)=
\sum_{n \in \Z} q^{n^2} z^{2n}
=1 + 2 \sum_{n=1}^{\infty} e^{\tau \pi i n^2} \cos (2 n \pi v).
\label{eqn:theta}
\end{align}
(Note that the present functions 
$\vartheta_{\mu}(v; \tau), \mu=1,2,3$ are denoted by
$\vartheta_{\mu}(\pi v,q)$,
and $\vartheta_0(v;\tau)$ by $\vartheta_4(\pi v,q)$ in \cite{WW27}.)
For $\Im \tau >0$, $\vartheta_{\mu}(v; \tau)$, $\mu=0,1,2,3$
are holomorphic for $|v| < \infty$
and satisfy the partial differential equation
\[
\frac{\partial \vartheta_{\mu}(v; \tau)}{\partial \tau}=
\frac{1}{4 \pi i} \frac{\partial^2 \vartheta_{\mu}(v; \tau)}{\partial v^2}.
\]
The parity with respect to $v$ is given by
\begin{equation}
\vartheta_1(-v; \tau)=-\vartheta_1(v; \tau),
\quad
\vartheta_{\mu}(-v; \tau)=\vartheta_{\mu}(v; \tau), \quad
\mu=0,2,3,
\label{eqn:even_odd}
\end{equation}
and they have the quasi-periodicity;  for instance, $\vartheta_1$ satisfies
\begin{equation}
\vartheta_1(v+1; \tau)=-\vartheta_1(v; \tau), \quad
\vartheta_1(v+\tau; \tau)=
-e^{-\pi i (2v+\tau)} \vartheta_1(v; \tau).
\label{eqn:quasi_periodic}
\end{equation}
By the definition (\ref{eqn:theta}), when $\tau \in \H$,
\begin{align}
& \vartheta_1(0; \tau)=\vartheta_1(1; \tau)=0, \qquad
\vartheta_1(x; \tau) > 0, \quad x \in (0,1),
\nonumber\\
& \vartheta_2(-1/2; \tau)=\vartheta_2(1/2; \tau)=0, \qquad
\vartheta_2(x; \tau) > 0, \quad x \in (-1/2, 1/2),
\nonumber\\
& \vartheta_0(x; \tau) > 0, \quad \vartheta_3(x; \tau) > 0, \quad x \in \R.
\label{eqn:theta_positive}
\end{align}
We see the asymptotics
\begin{align}
& \vartheta_0(v; \tau) \sim 1, \quad
\vartheta_1(v; \tau) \sim 2 e^{\tau \pi i/4} \sin (\pi v), \quad
\vartheta_2(v; \tau) \sim 2 e^{\tau \pi i/4} \cos(\pi v), \quad
\vartheta_3(v; \tau) \sim 1,
\nonumber\\
&\qquad \qquad \qquad \mbox{in} \quad
\Im \tau \to + \infty \quad
({\it i.e.}, \quad q=e^{\tau \pi i} \to 0).
\label{eqn:theta_asym}
\end{align}
The following functional equalities are known as
Jacobi's imaginary transformations \cite{WW27,NIST10},
\begin{align}
\vartheta_0(v; \tau)
&= e^{\pi i/4} \tau^{-1/2} e^{-\pi i v^2/\tau}
\vartheta_2 \left( \frac{v}{\tau}; - \frac{1}{\tau} \right),
\nonumber\\
\vartheta_1(v; \tau)
&= e^{3 \pi i/4} \tau^{-1/2} e^{-\pi i v^2/\tau}
\vartheta_1 \left( \frac{v}{\tau}; - \frac{1}{\tau} \right),
\nonumber\\
\vartheta_2(v; \tau)
&= e^{\pi i/4} \tau^{-1/2} e^{-\pi i v^2/\tau}
\vartheta_0 \left( \frac{v}{\tau}; - \frac{1}{\tau} \right), 
\nonumber\\
\vartheta_3(v; \tau)
&= e^{\pi i/4} \tau^{-1/2} e^{-\pi i v^2/\tau}
\vartheta_3 \left( \frac{v}{\tau}; - \frac{1}{\tau} \right). 
\label{eqn:Jacobi_imaginary}
\end{align}

\SSC{Selberg-type Integral Formulas
Including the Jacobi Theta Functions}
\label{sec:appendixB}
Apply the Macdonald denominator formulas
(\ref{eqn:Macdonald3}) with (\ref{eqn:a1}) into (\ref{eqn:q1}).
Then Lemma \ref{thm:normalization} gives the following
Selberg-type integral formulas
including the Jacobi theta functions.
Let $s(N)=0$ if $N$ is even, and $s(N)=3$
if $N$ is odd.
For $0 < t_{\ast} < \infty$, $t \in (0, t_{\ast})$,
\begin{align*}
&
\int_{[0, 2\pi r]^N} d \x \,
\vartheta_{s(N)} \left( \sum_{j=1}^N \xi(x_j) ;
\cN^{\AN} \tau(t_{\ast}-t) \right)
\vartheta_{s(N)} \left( \sum_{j=1}^N \xi(x_j);
\cN^{\AN} \tau(t) \right) 
\nonumber\\
& \qquad 
\times
W^{\AN}(\xi(\x); \cN^{\AN} \tau(t_{\ast}-t))
W^{\AN}(\xi(\x); \cN^{\AN} \tau(t))
= \frac{N! \prod_{n=1}^N m^{\AN}_n(t_{\ast})}{a^{\AN}(t_{\ast}-t) a^{\AN}(t)},
\nonumber
\\
& 
\int_{[0, \pi r]^N} d \x \,
W^{\RN}(\xi(\x); \cN^{\RN} \tau(t_{\ast}-t))
W^{\RN}(\xi(\x); \cN^{\RN} \tau(t))
= \frac{N! \prod_{n=1}^N m^{\RN}_n(t_{\ast})}{a^{\RN}(t_{\ast}-t) a^{\RN}(t)},
\nonumber\\
& \qquad \qquad \qquad \qquad
\mbox{for $\RN=\BN,\BNv,\CN,\CNv,\BCN,\DN$}.
\label{eqn:SelbergBD}
\end{align*}
In particular, if we set $t=t_{\ast}/2$, we have the following,
\begin{align*}
&
\int_{[0, 2\pi r]^N} d \x \,
\left\{ 
\vartheta_{s(N)} \left( \sum_{j=1}^N \xi(x_j) ;
\cN^{\AN} \tau(t_{\ast}/2) \right)
W^{\AN}(\xi(\x); \cN^{\AN} \tau(t_{\ast}/2))
\right\}^2
\nonumber\\
& \qquad \qquad \qquad \qquad
= \frac{N! \prod_{n=1}^N m^{\AN}_n(t_{\ast})}{ (a^{\AN}(t_{\ast}/2))^2},
\nonumber
\\
& 
\int_{[0, \pi r]^N} d \x \,
\left\{ W^{\RN}(\xi(\x); \cN^{\RN} \tau(t_{\ast}/2)) \right\}^2
= \frac{N! \prod_{n=1}^N m^{\RN}_n(t_{\ast})}{ (a^{\RN}(t_{\ast}/2))^2},
\nonumber\\
& \qquad \qquad \qquad \qquad
\mbox{for $\RN=\BN,\BNv,\CN,\CNv,\BCN,\DN$}.
\end{align*}

\SSC{Determinantal Point Processes
and Correlation Kernels}
\label{sec:appendixC}

Let $N \in \N$, $S \subset \R^d$.
Assume that the probability measure of point process, 
$\Xi(\cdot)=\sum_{j=1}^N \delta_{X_j}(\cdot)$, is given by
\begin{equation}
\bP(\X \in d \x) = \bp(\x) d \x
=\frac{1}{C(N)} \det_{1 \leq j, k \leq N} [f_j(x_k)] 
\det_{1 \leq \ell, m \leq N} [g_{\ell}(x_m)],
\quad \x \in S^N,
\label{eqn:P_det1}
\end{equation}
with the biorthogonality relations
\begin{equation}
\int_S f_j(x) g_k(x) dx
= h_j \delta_{jk}, \quad j, k \in \{1,2, \dots, N\},
\label{eqn:ortho1}
\end{equation}
where $h_j > 0, j=1,2, \dots, N$.
Let $\rC_{\rc}(S)$ be a collection of all continuous real
functions with a compact support in $S$.
For $\psi \in \rC_{\rc}(S)$, $\theta \in \R$,
the characteristic function of $(\Xi, \bP)$ is defined as
\[
\Psi[\psi; \theta]
= \bE\left[ e^{\theta \sum_{j=1}^N \psi(X_j)} \right]
= \frac{1}{N!} \int_{S^N} d \x \,
e^{\theta \sum_{j=1}^N \psi(x_j)} \bp(\x),
\]
which can be regarded as the Laplace transform of $\bp$.
Put $\chi(x)=1-e^{\theta \psi(x)}$.
By performing binomial expansion, we obtain
\begin{align}
\Psi[\psi; \theta]
&= \frac{1}{N!} \int_{S^N} d \x \,
\prod_{j=1}^N (1-\chi(x_j)) \bp(\x)
\nonumber\\
&= 1+\sum_{n=1}^N (-1)^n \frac{1}{n!}
\int_{S^n} \prod_{\ell=1}^n \{ dx_{\ell} \chi(x_{\ell}) \}
\rho(\{x_1, \dots, x_n \}),
\label{eqn:Psi1}
\end{align}
where $\rho(\{x_1, \dots, x_n\})$ is given by (\ref{eqn:correlation1}).
This implies that if we regard $\Psi[\psi; \theta]$ as
a functional of $\chi$, it gives the generating function of
correlation functions \cite{KT07}.
Insert (\ref{eqn:P_det1}) into (\ref{eqn:Psi1}) and use the
Heine identity
\begin{equation}
\frac{1}{N!} \int_{S^N} d \x \,
\det_{1 \leq j, k \leq N} [ \phi_j(x_k) ]
\det_{1 \leq \ell, m \leq N} [ \varphi_{\ell}(x_m)]
=\det_{1 \leq j, k \leq N} 
\left[ \int_{S} d x \, \phi_j(x) \varphi_k(x) \right]
\label{eqn:Heine}
\end{equation}
for square integrable functions $\phi_j, \varphi_j, j \in \{1,2, \dots, N\}$.
Then we have
\begin{align*}
\Psi[\psi; \theta]
&= \frac{1}{C(N)} 
\det_{1 \leq j, k \leq N} \left[
\int_{S} dx \, f_j(x)(1-\chi(x)) g_k(x) \right]
\nonumber\\
&= \frac{\det_{1 \leq j, k \leq N} \left[
\int_{S} dx \, f_j(x)(1-\chi(x)) g_k(x) \right]}
{\det_{1 \leq j, k \leq N} \left[
\int_{S} dx \, f_j(x)g_k(x) \right]}
\nonumber\\
&= \frac{\det_{1 \leq j, k \leq N} \left[
\int_{S} dx \, f_j(x) g_k(x) 
- \int_{S} dx \, f_j(x) \chi(x) g_k(x)  \right]}
{\det_{1 \leq j, k \leq N} \left[
\int_{S} dx \, f_j(x)g_k(x) \right]},
\end{align*}
where we used the normalization condition
$\Psi[\psi; 0]=1$ at the second equality.
We introduce the $N \times N$ matrices $\sfA$ and
$\sfA[\chi]$ with the entries
\[
(\sfA)_{jk}=\int_{S} dx \, f_j(x) g_k(x), \quad
(\sfA[\chi])_{jk} = \int_{S} dx \, f_j(x) \chi(x) g_k(x),
\quad 1 \leq j, k \leq N. 
\]
Since (\ref{eqn:ortho1}) is assumed, $\sfA$ is a regular matrix,
and the above is written as
\[
\Psi[\psi; \theta]
=\det_{1 \leq j, k \leq N} 
\left[ \delta_{jk} - (\sfA^{-1} \sfA[\chi])_{jk} \right], 
\]
where
\[
(\sfA^{-1} \sfA[\chi])_{jk}= \int_{S} dx B_j(x) g_k(x)
\quad 
\mbox{with}
\quad
B_j(x)=\sum_{\ell=1}^N (\sfA^{-1})_{j \ell} f_{\ell}(x) \chi(x).
\]
Now we apply the Fredholm expansion formula.
Then we can verify that \cite{KT07}
\[
\Psi[\psi; \theta]
= 1 + \sum_{n=1}^N (-1)^n
\frac{1}{n!} \int_{S^n} \prod_{\ell=1}^n dx_{\ell} 
\det_{1 \leq j, k \leq n} 
\left[ \sum_{m=1}^N g_m(x_j) B_m(x_k) \right].
\]
By the orthogonality (\ref{eqn:ortho1}), 
$(A^{-1})_{j \ell}=(1/h_j) \delta_{j \ell}$ and hence
$B_m(x)=f_m(x) \chi(x)/h_m$.
We define
\begin{equation}
K(x,y)=\sum_{m=1}^N \frac{g_m(x) f_m(y)}{h_m},
\quad x, y \in S.
\label{eqn:K1}
\end{equation}
Then we arrive at the expression
\begin{equation}
\Psi[\psi; \theta]
= 1+ \sum_{n=1}^N (-1)^n \frac{1}{n!}
\int_{S^n} \prod_{\ell=1}^n \{ dx_{\ell} \chi(x_{\ell}) \}
\det_{1 \leq j, k \leq n} [K(x_j, x_k)].
\label{eqn:Psi2}
\end{equation}
For any $\chi \in \rC_{\rc}(S)$, we have proved that 
(\ref{eqn:Psi1}) is equal to
(\ref{eqn:Psi2}). Hence we can conclude that
\[
\rho(\{x_1, \dots, x_n \})
=\det_{1 \leq j, k \leq n} [ K(x_j, x_k)],
\quad n = 1,2, \dots, N.
\]
In summary, if (\ref{eqn:P_det1}) and (\ref{eqn:ortho1}) are satisfied,
then the point process $(\Xi, \bP)$ is determinantal
and the correlation kernel is given by (\ref{eqn:K1}).
We note that (\ref{eqn:Psi2}) defines the
{\it Fredholm determinant} associated with the integral kernel
$K(x, y) \chi(y)$  
and it is written as
$\Det_{x, y \in S} \Big[\delta(x-y)-K(x,y) \chi(y) \Big]$ 
(see, for instance, \cite{Kat15_Springer}).



\end{document}